\documentclass[11pt]{article}

\usepackage{amssymb,graphics,amsmath,amsthm,amsopn,amstext,amsfonts,color,fullpage,hyperref,bm}
\usepackage{subfig,placeins} % make it possible to include more than one captioned figure/table in a single float
\usepackage{graphicx}

\newcommand{\gap}{\vspace{0.1in}}
\newcommand{\epc}{\hspace{1pc}}
\newcommand{\thalf}{{\textstyle{\frac{1}{2}}}}
\newcommand{\wh}{\widehat}
\newcommand{\wt}{\widetilde}

\newcommand{\vect}[1]{\mathbf{#1}}

\newtheorem{theorem}{Theorem}
\newtheorem{proposition}[theorem]{Proposition}

\newtheorem{corollary}[theorem]{Corollary}
\newtheorem{definition}[theorem]{Definition}
\newtheorem{example}[theorem]{Example}

\newcounter{algo}

\title{Computing B-Stationary Points of Nonsmooth DC 
Programs\footnote{This work was based on research partially supported by the
U.S.\ National Science Foundation grants CMMI 1402052 and 0969600.}}

\date{Original: October 2014; revised August 2015}

\author{Jong-Shi Pang\footnote{Department of Industrial and Systems Engineering, University of Southern
California, Los Angeles, California  90089-0193, U.S.A.
{\sl Email:} jongship@usc.edu.} \and
Meisam Razaviyayn\footnote{Department of Electrical and Computer Engineering,
University of Minnesota, Minneapolis, Minnesota 55455, U.S.A.
{\sl Email:} razav002@umn.edu} \and Alberth Alvarado\footnote{Department of Industrial and Enterprise Systems Engineering,
University of Illinois at Urbana-Champaign, Urbana, Illinois 61801, U.S.A.
{\sl Present address: Department of Applied Mathematics, Universidad Galileo, Guatemala, 01010, Guatemala}.
{\sl Email:} alberth@galileo.edu}}

\begin{document}

\maketitle

\begin{abstract}
\noindent Motivated by a class of applied problems arising from physical layer based security in a digital communication
system, in particular, by a secrecy sum-rate maximization problem, this paper studies a nonsmooth, difference-of-convex (dc)
minimization problem.
The contributions of this paper are: (i) clarify several kinds of stationary solutions and their relations;
(ii) develop and establish the convergence of a novel algorithm for computing a d-stationary solution of a problem
with a convex feasible set that is arguably the sharpest kind among the various stationary solutions; (iii) extend
the algorithm in several directions including: a randomized choice of the subproblems that could help the
practical convergence of the algorithm, a distributed penalty approach for problems whose objective functions are sums
of dc functions, and problems with a specially structured (nonconvex) dc constraint.  For the latter class of problems,
a pointwise Slater
constraint qualification is introduced that facilitates the verification and computation of a B(ouligand)-stationary point.
\end{abstract}

\section{Introduction}

A general difference-of-convex (dc) optimization problem refers to the minimization/maximization of an objective function
that is the difference of two convex functions subject to constraints defined by functions
of the same kind.  Such optimization problems form a large class of nonconvex programs and have been
studied extensively for more than three decades in the mathematical programming
literature \cite{ATao05,Hiriart-Urruty85,HThoai99,Moudafi08,MMoudafi08,MMoudafi06,SLanckriet12,TTAn97,Tuy87}.
In particular, Pham Thinh Tao and Le Thi Hoai An have made pioneering contributions to this important subfield
of contemporary optimization and are responsible for much of the development of theory, algorithms, and applications
of dc programming to date.  See the cited references \cite{ATao05,HoiAnTao14,HoiAnTao14-1,HThoai99,HoiAnTaoHuynhi12,TTAn97} for a
sample of their voluminous writings in this area.
In particular, the DCA (Difference-of-Convex Algorithm) has been a principal algorithm for computing a {\sl critical point}
of the problem.

\gap

Our interest in this class of nonconvex optimization problems stemmed initially from a particular application
pertaining to physical layer based security in a
digital communication system \cite{Alvarado14,ASPang14} and a related one of joint base-station assignment and
power allocation \cite{SRLuo14}.  A first glance at their formulations
% involving a bivariate function
does not immediately reveal that the resulting nonsmooth maximization problem (see (\ref{eq:general_problem})) is of the dc
type.  Yet, a careful look at the objective function
shows that it can be expressed as the difference of two concave functions, one of which is differentiable and the other
one is not.  Furthermore, via a ``lifting'' of the problem using some auxiliary continuous variables to
express a discrete pointwise maximum function as a value function of an optimization problem over the unit simplex,
this applied problem can be formulated as a smooth, bi-concave (thus nonconcave), linearly constrained maximization problem.
This special problem raises several interesting questions that do not seem to have been adequately addressed in the existing
literature of dc programming.  As a linearly constrained dc program, one can speak about the concept of a d(irectional)-stationary
point of the problem, i.e., a point at which the one-sided derivatives of the objective function along any feasible
(equivalently, tangent) direction are
nonnegative.  Since the lifted formulation is smooth, one can speak about the standard concept of stationarity, which we call
{\sl lifted stationarity}, in terms of the gradient of the objective function in the lifted space.  The following is a set
of questions that have partially motivated our research: for a dc program with convex constraints,
\begin{description}
\item[\rm (a)] How are the concepts of a critical point and d-stationarity related to each other?
\item[\rm (b)] How is lifted stationarity defined in general?  How is it related to criticality and d-stationarity of the un-lifted
problem?
\item[\rm(c)] Are there algorithms that can provably compute a d-stationary point?
\end{description}
Providing answers to the first two questions constitute a major part of our study.  In so doing, we are led to the contention that d-stationarity
is arguably the sharpest among these stationarity concepts and yet the computation of such a point by an existing provably convergent
algorithm seems to have been elusive to date.
This lack of a computational scheme for obtaining a d-stationary point of a convex-constrained dc program leads to the other part of
our work, namely, to propose a novel iterative algorithm to fill this
gap.  The design of the algorithm is interesting in its own right, namely, it contains innovative ideas that do not seem to have been
introduced in the dc literature; in particular, we present a randomized version of the algorithm to deal with a potential weakness
of the algorithm in practical implementation.  Convergence of the algorithms is established.

\gap

Also included in our algorithmic development is the extension of the basic algorithm to a multi-agent context where the objective is
the sum of dc functions with each summand being a private objective function (with coupled variables) of an individual agent.
[The applied problem mentioned at
the beginning of this introduction is a problem of this kind.]  In such a context, it is desirable to develop a distributed algorithm
wherein the optimization of each agent can be carried out independently of the other agents.  The design of such a distributed algorithm
is another major contribution of our work. This is accomplished via a double iteration wherein the outer loop is a penalty-based scheme
and each inner loop applies the newly developed algorithm for computing a d-stationary point of a penalized subproblem.  The separability
into individual agent-based optimization occurs naturally in the latter loop.

\gap

Our last contribution is the extension of the basic algorithm
to allow for the presence of a non-differentiable dc constraint, leading what has been called a {\sl general dc program}
\cite{HoiAnTao14-1,HoiAnHuynhTao14,TTAn97}.  Such a constraint adds considerable complication to the theory
and computation for a convex feasible set, due to the nonsmoothness and nonconvexity of the dc constraint.  For such a dc constrained
dc program, we formally define the concept of a B(ouligand)-stationary point
and show how it can be characterized by a reasonable number of convex programs, thus making the verification of such a stationary
point practically implementable.  A provably convergent algorithm is then developed for computing such a point.
% Application of the dc constrained dc program to optimization problems with complementarity constraints will be
% briefly sketched, leaving the details for further study and to be reported in a separate work.

\section{Motivating Applied Problems}

In this section, we discuss two applied problems pertaining to power allocation in digital communication systems
that had motivated our research.
These problems lead us to a unified class of value functions of a continuum family of bivariate functions which we show are of the dc type.
For more applications of dc programs to communication
systems and other domains, we refer the readers to \cite{HoiAnTao14,ATao05}.  Subsequent to the completion of
the original version of this paper, the authors recognized that many interesting classes of nonconvex optimization problems
are actually of the dc kind and have not been treated as such; examples include those arising from deviation
measures in risk analysis as well as in the regularization of loss functions in statistical
learning.  Due to space limitations, we cannot give details of these other problems.  We are hopeful that our work herein
opens renewed opportunities for the dc methodology, in particular the results and algorithms in this paper and those existed
in the literature, to be applied to deal with these nonconvex problems more effectively.
% We plan to investigate such opportunities in detail and will report our findings elsewhere.

% We will derive some stationarity results for these problems that will subsequently be generalized to an abstract setting.
% We believe that the extended results will have broad applications to dc programs arising from other sources.

\gap

The concept of secrecy capacity is of fundamental importance in information theory \cite{JWGerbracht10}.
Based on this concept, a design problem in physical layer based security is to allocate power budget to the network spectra
so that the transmissions between legitimate parties can be kept secure.  The problem stated below pertains to the
single-input-single-output (SISO) paradigm
where users of the network (consisting of
transmitter-receiver pairs) communicate over multiple non-orthogonal subchannels.  There are also a number of ``friendly''
jammers and one eavesdropper.   Each legitimate user's transmitter wants to communicate (in a secure way) with
its corresponding receiver over a set of parallel subchannels. The friendly jammers are entities willing to
cooperate with the legitimate parties by introducing judicious interferences so as to impair the eavesdropper's ability
to decode the messages between intended nodes.  With $H_{rq}(k)$, $H_{r0}(k)$ and $\wh{H}_{j0}(k)$ denoting
the channel gains and $\sigma_q^2(k)$ the variances of channel noise, all being constants in the model, and
$p_q(k)$ and $\wh{p}_j(k)$ denoting, respectively, the variable power of user $q$ and jammer $j$
allocated to channel $k$, this multi-jammer secrecy-sum-rate maximization problem is formulated as follows
\cite{ASPang14}:
\begin{equation}\label{eq:SRsingleProblem}
\begin{array}{ll}
\displaystyle{
{\operatornamewithlimits{\mbox{maximize}}_{\left( \mathbf{p}, \mathbf{\wh{p}} \right) \, \geq \, \mathbf{0}}}
} & \displaystyle{
\sum_{q=1}^Q
} \, \displaystyle{
\sum_{k=1}^N
} \, \left[ \, R_{qqk}(\mathbf{p},\mathbf{\wh{p}}) - R_{q0k}(\mathbf{p},\mathbf{\wh{p}}) \, \right]^+ \\ [0.2in]
\mbox{subject to:} & \displaystyle{
\sum_{k=1}^{N}
} \, p_q(k) \, \leq \,  P_q^{\max}  \,\,\,\,\,\,\,\, \forall \, q=1,\ldots,Q \,\,\,\,\,\,\,\ \mbox{(agents' private constraints) } \\ [0.2in]
\mbox{and} & \displaystyle{
\sum_{k=1}^{N}
} \, \wh{p}_j(k) \, \leq \, \wh{P}_j^{\max} \,\,\,\,\,\,\,\, \forall \,  j=1,\ldots,J \,\,\,\,\,\,\,\ \mbox{ (coupling constraints) },
\end{array}
\end{equation}
where $\vect p \triangleq \left( \left( p_q(k) \right)_{k=1}^{N} \right)_{q=1}^{Q}$,
$\wh{\vect p} \triangleq \left( \left( \wh{p}_j(k) \right)_{k=1}^{N} \right)_{j=1}^{J}$,
$[ \, \bullet \, ]^+ \triangleq \max(0,\bullet)$ is the plus-function, and
$R_{qqk}(\mathbf{p},\mathbf{\wh{p}})$ and $R_{q0k}(\mathbf{p},\mathbf{\wh{p}})$ are
the Shannon information rate functions given by:
\[ \begin{array}{lll}
R_{qqk}(\mathbf{p},\mathbf{\wh{p}}) & \triangleq & \log \left( 1+ \displaystyle{
\frac{H_{qq}(k) \, p_q(k)}{\sigma_q^2(k) + \displaystyle{
\sum_{q \neq r =1}^Q
} \, H_{rq}(k) \, p_r(k) + \displaystyle{
\sum_{j =1}^{J}
} \, \wh{H}_{jq}(k) \, \wh{p}_j(k)}
} \, \right) \\ [0.6in]
R_{q0k}(\mathbf{p},\mathbf{\wh{p}}) & \triangleq & \log \left( 1+ \displaystyle{
\frac{H_{q0}(k) \, p_q(k)}{\sigma_q^2(k) + \displaystyle{
\sum_{q \neq r=1}^Q
} \, H_{r0}(k) \, p_r(k) + \displaystyle{
\sum_{j =1}^{J}
} \, \wh{H}_{j0}(k) \, \wh{p}_j(k)}
} \, \right);
\end{array} \]
both $R_{qqk}$ and $R_{q0k}$ are clearly differentiable differences of concave functions.
Since the pointwise maximum of a finite number of dc functions is a dc function
\cite{HThoai99,TTAn97}, and so is the sum of a finite number of dc functions, it follows
that the objective function of (\ref{eq:SRsingleProblem}) is a dc function.

\gap

A related problem is that of optimal joint base station assignment
and power allocation in a communication network \cite{SRLuo14}.
Admitting a similar formulation with binary variables subject to
a knapsack constraint, thus with multiple (more than two) discrete choices,
this problem is 
\begin{equation} \label{eq:jointassall}
\begin{array}{ll}
\displaystyle{
\operatornamewithlimits{\mbox{maximize}}_{x_q,y_q}
} & \displaystyle{
\sum_{q=1}^Q
} \, \displaystyle{
\sum_{\ell=1}^L
} \, y_{q\ell} \, \left[ \, \displaystyle{
\sum_{k=1}^N
} \, \log\left( \, 1 + \displaystyle{
\frac{H_{q\ell}(k) \, x_{q\ell}(k)}{\sigma_\ell(k)^2 + \displaystyle{
\sum_{r \neq q}
} \, y_{r\ell} \, H_{r \ell}(k) \, x_{rk}(f)}
} \, \right) \, - \underbrace{c_{q\ell}}_{\mbox{set-up cost}} \, \right] \\ [0.4in]
\mbox{subject to} & \displaystyle{
\sum_{\ell=1}^L
} \, \displaystyle{
\sum_{k=1}^N
} \, x_{q\ell}(k) \, \leq \, B_q^{\max} \\ [0.3in]
\mbox{and} & \mbox{for all $q \, = \, 1, \cdots, Q$}, \\ [0.1in]
& 0 \, \leq \, x_{q\ell}(k) \, \leq \, \mbox{CAP}_{q\ell}(k), \epc \forall \, \ell \, = \, 1, \cdots, L
\mbox{ and } k \, = \, 1, \cdots, N \\ [5pt]
& \displaystyle{
\sum_{\ell=1}^L
} \, y_{q\ell} \, = \, 1, \epc y_{q\ell} \, \in \, \{ \, 0,1 \, \}, \epc \forall \, \ell \, = \, 1, \cdots, L,
\end{array}
\end{equation}
where the additional index $\ell = 1, \cdots, L$ labels the base stations.  The problem
is equivalent to
\[ \begin{array}{ll}
\displaystyle{
\operatornamewithlimits{\mbox{maximize}}_{x_q}
} & \theta(x) \, \triangleq \, \displaystyle{
\sum_{q=1}^Q
} \, \underbrace{\displaystyle{
\operatornamewithlimits{\mbox{maximum}}_{y_q \in Y_q}
} \, \displaystyle{
\sum_{\ell=1}^L
} \, y_{q\ell} \, \left[ \, \displaystyle{
\sum_{k=1}^N
} \, \log\left( \, 1 + \displaystyle{
\frac{H_{q\ell}(k) \, x_{q\ell}(k)}{\sigma_\ell(k)^2 + \displaystyle{
\sum_{r \neq q}
} \, y_{r\ell} \, H_{r \ell}(k) \, x_{rk}(f)}
} \, \right) \, - c_{q\ell} \, \right]}_{\mbox{pointwise max of dc functions}} \\ [0.5in]
\mbox{subject to} & \mbox{for all $q = 1, \cdots, Q$}, \\ [0.1in]
& \displaystyle{
\sum_{\ell=1}^L
} \, \displaystyle{
\sum_{k=1}^N
} \, x_{q\ell}(k) \, \leq \, B_q^{\max} \\ [0.3in]
\mbox{and} &
0 \, \leq \, x_{q\ell}(k) \, \leq \, \mbox{CAP}_{q\ell}(k), \epc \forall \, \ell \, = \, 1, \cdots, L
\mbox{ and } k \, = \, 1, \cdots, N,
\end{array} \]
where $Y_q \triangleq \left\{ y_q \in \{ 0,1 \}^L \mid \displaystyle{
\sum_{\ell=1}^L
} \, y_{q\ell} = 1 \right\}$ remains a discrete set that contains a sum constraint which
in some applications could be generalized to a cardinality constraint of the type:
$K \geq \displaystyle{
\sum_{\ell=1}^L
} \, y_{q\ell} \geq 1$.  Since each $Y_q$ is a discrete set, it follows again that the above objective
$\theta(x)$ is a dc function.

\subsection{A digression: Continuum family of dc functions}

While it is known that the pointwise maximum of finitely many dc functions is a dc function, it is not
known whether this dc property extends to the pointwise maximum of a continuum family of dc
functions.  It turns out that this extension does not hold in general as the example below
shows.

\begin{example} \rm
Let $g : \mathbb{R}^n \to \mathbb{R}$ be a (globally) Lipschitz continuous function that is not directionally
differentiable everywhere.  Such a function exists as provided by a component function of the Euclidean projector
onto a specially constructed compact convex set (see \cite{Shapiro94}
for $n = 2$ and \cite{Kruskal69} and \cite[Exercise~4.8.5]{FPang03} for $n = 3$).  Let $L > 0$ be a Lipschitz
constant of $g$ in the $\ell_1$-norm; i.e.,
\[
| \, g(x) - g(y) \, | \, \leq \, L \, \| \, x - y \, \|_1, \epc \forall \, x, \, y \, \in \, \mathbb{R}^n.
\]
Define $f(x,y) \triangleq g(y) - L \, \| x - y \|_1$.  Obviously, $g(x) = \displaystyle{
\operatornamewithlimits{\mbox{maximum}}_{y \in \mathbb{R}^n}
} \, f(x,y)$.  It is clear that $f(\bullet,y)$ is a concave, thus dc, function.  Yet $g$ cannot be a dc function
as every dc function must be directionally differentiable but $g$ is not.  \hfill $\Box$
\end{example}

In what follows, we present a class of value functions of bivariate functions that preserves the dc property;
it turns out the structure of the component functions is important; such a structure includes the case of
finitely many dc functions, and in particular the two problems (\ref{eq:SRsingleProblem}) and (\ref{eq:jointassall}).
Specifically, consider the following non-convex, non-differentiable multi-agent optimization problem:
\begin{equation}\label{eq:general_problem}
	\underset{x \in X}{\text{ maximize}} \ \theta(x) \, \triangleq \,
		\sum_{i=1}^I  \left(
		\underset{\lambda^i \in \Lambda^i}{\text{maximum}} \,
		\sum_{j=1}^{J} h_{i,j} (\lambda^i) f_{i,j}(x)
		\right),
	%\right\},
\end{equation}
where the set $X \subseteq \mathbb{R}^{n}$ is closed and convex, and for each
$i=1,\cdots,I$ (denoting the agents' labels), $\Lambda^i \subseteq \mathbb{R}^{m_i}$
is a compact set (not assumed to be convex; cf.\ e.g., the set $Y_q$ in the previous
subsection).  For each $i=1,\cdots,I$ and $j=1,\cdots,J$,
$f_{i,j} : \Omega \subseteq \mathbb{R}^n \rightarrow \mathbb{R}$, where $\Omega$ is an
open convex superset of $X$, is  
% a continuously differentiable function defined on the open convex superset $\Omega$ of $X$.
% Moreover, we assume that each $f_{i,j}$ is 
either convex or concave on $X$.
Finally, for each $i=1,\cdots,I$ and $j=1,\cdots,J$, $h_{i,j} : \Omega^{\, i} \subseteq \mathbb{R}^{m_i} \rightarrow \mathbb{R}$,
where $\Omega^{\, i}$ is open convex set containing $\Lambda^i$, is such that each product
$h_{i,j} (\lambda^i) f_{i,j}(x)$ is concave in $\lambda^i$ for fixed $\vect x$.
A particularly important special case of $\Lambda^i$ is when it is a unit simplex and each function $h_{i,j}$ is affine
so that the continuous pointwise maximum becomes a discrete pointwise maximum and the overall problem (\ref{eq:general_problem})
is as follows:
\[
\displaystyle{
\operatornamewithlimits{\mbox{maximize}}_{x \in X}
} \ \displaystyle{
\sum_{i=1}^I
} \, \displaystyle{
\max_{1 \leq j \leq J}
} \, f_{i,j}(x).
\]

\gap

Our proof showing that the function $\theta$ in (\ref{eq:general_problem}) is
of the dc kind appears to be new.  In order not to further deviate from the discussion of the main topics of
this paper, we provide the proof in an appendix at the end of the paper.  Notice that (\ref{eq:general_problem})
is ``equivalent'' to the ``bivariate'' maximization
\begin{equation}\label{eq:restated general_problem}
	\underset{x \in X; \, ( \lambda^i \in \Lambda^i )_{i=1}^I}{\text{ maximize}} \
\displaystyle{
\sum_{i=1}^I
} \, \displaystyle{
\sum_{j=1}^{J}
} \, h_{i,j} (\lambda^i) f_{i,j}(x),
\end{equation}
where the equivalence pertains to the globally optimal solutions of these two problems.  Nevertheless, when it comes to
stationary solutions, the situation is quite different; see the subsequent discussion in Subsection~\ref{subsec:lifted}.
In particular, (\ref{eq:restated general_problem}) is a differentiable program if all function $h_{i,j}$ and $f_{i,j}$
are differentiable while (\ref{eq:general_problem}) remains non-differentiable due to the max operator; even with this advantage of
differentiability, the concept of ``d-stationarity'' in the former problem is not as sharp as the same concept in the latter that
has the $\lambda$-variable ``hidden'', i.e., in the $x$-alone formulation:
\[
\underset{x \in X}{\text{ maximize}} \ \theta(x) \, \triangleq \,
		\sum_{i=1}^I \, \theta_i(x); \epc \mbox{with each} \epc
\theta_i(x) \, \triangleq \, \underset{\lambda^i \in \Lambda^i}{\text{maximum}} \,
		\sum_{j=1}^{J} h_{i,j} (\lambda^i) f_{i,j}(x).
\]
The upshot of this discussion is that different formulations of a non-differentiable, non-convex optimization problem could lead
to stationary solutions with unequal likelihood for being a locally optimal solution.  The search for a superior formulation
is not an easy task in general, however.

% Even this result relies critically on the structure of each summand
% within the maximand
% in $\theta_i$, i.e., the product of two univariate functions possessing further properties as specified above.

\section{Stationarity: Convex Constraints}

As a non-convex optimization problems, globally optimal solutions of a dc program are in general not possible
to be computed.  Thus, one has to settle for computing a ``stationary'' solution in practice.  Even so, one
has to be cautious about the notion of stationarity, especially in the case where the constraints contain dc
functions.  The situation is simpler when the constraint set is convex; in this section, we consider this case
first.   Specifically, we deal with the following convex constrained dc minimization program:
\begin{equation} \label{eq:ccdc}
\displaystyle{
{\operatornamewithlimits{\mbox{minimize}}_{x \in X}}
} \ \zeta(x) \, \triangleq \, f(x) - g(x),
\end{equation}
where $f$ and $g$ are convex functions defined on an open convex set $\Omega$ containing the closed convex set $X \subseteq \mathbb{R}^n$.
[Note the change from maximization in the previous section to minimization in the problem (\ref{eq:ccdc}).]
Since $\zeta$ is not differentiable, stationarity concepts of (\ref{eq:ccdc}) are defined in terms of directional
derivatives of the objective function, which we briefly review in the subsection below.  Before doing so, we
mention the references \cite{BBorwein11,Hiriart-Urruty85} where a host of properties of dc functions are summarized.

\subsection{Directional derivatives}

The directional derivative of $\zeta$ at a point $x \in \Omega$ along a direction $d \in \mathbb{R}^n$ is given by:
\[
\zeta^{\, \prime}(x;d) \, \triangleq \, \displaystyle{
\lim_{\tau \downarrow 0}
} \, \displaystyle{
\frac{\zeta(x + \tau d) - \zeta(x)}{\tau}
} .
\]
It is well known that convex functions are directionally differentiable; so with $\zeta = f - g$ being a dc function,
$\zeta^{\, \prime}(x;d)$ is well defined for all $x \in \Omega$ and $d \in \mathbb{R}^n$; moreover
\[
\zeta^{\, \prime}(x;d) \, = \, f^{\, \prime}(x;d) - g^{\, \prime}(x;d).
\]
Since dc functions are locally Lipschitz continuous (and are thus B(ouligand) differentiable \cite[Definition~3.1.2]{FPang03}),
the C(larke) directional derivative is also well defined:
\[
\zeta^0(x;d) \, \triangleq \, \displaystyle{
\limsup_{\stackrel{y \to x}{\tau \downarrow 0}}
} \, \displaystyle{
\frac{\zeta(y + \tau d) - \zeta(y)}{\tau}
} \, .
\]
In general, $\zeta^0(x;d) \geq \zeta^{\, \prime}(x;d)$.  These two directional derivatives are equal if the function
$\zeta$ is C-regular \cite{Clarke87}.  However, dc functions are in general not C-regular.  We recall that a function
$\zeta$ is {\sl strictly differentiable} at a point $x$ if the following limit holds:
\[
\displaystyle{
\lim_{\stackrel{(y,z) \to (x,x)}{y \neq z}}
} \, \displaystyle{
\frac{\zeta(y) - \zeta(z) - \nabla \zeta(x)^T ( \, y - z \, )}{\| \, y - z \, \|}
} \, = \, 0,
\]
where $\nabla \zeta(x)$ denotes the gradient vector of $\zeta$ at $x$.
If $g$ is strictly differentiable at $x$, then the dc function $\zeta$ is C-regular at $x$.  This class of dc functions
deserves a name.  Specifically, we say that $\zeta$ is a {\sl good} dc function on $\Omega$ (with respect to
a minimization problem) if there exists a strictly differentiable convex function $v$ on $\Omega$ such that $\zeta + v$ is
convex on $\Omega$; in other words, $\zeta$ is a good dc function if convex functions $u$ and $v$
exist such that $\zeta = u - v$ and $v$ is strictly differentiable.  The class of good dc functions appears
extensively in the machine learning area; see e.g.\ \cite{SLanckriet12} and the references therein.
These dc functions are particularly relevant in the context of computing stationary solutions and
play an important role in the convergence of several families of iterative algorithms for solving dc programs,
such as: the DCA \cite{ATao05,HThoai99,SLanckriet12,TTAn97} that has been a
fundamental algorithm with many applications,
the S(uccessive)C(onvex)A(approximation) method \cite{ASPang14,BSTeboulle13,Hong2013Linear,RHLPang14,SFLSong14,SFSPPang14}
that has
attracted significant interest in recent years for solving non-convex non-differentiable optimization problems,
and an alternating/successive minimization method \cite{RHLuo13,RHLPang14}
for solving the joint minimization formulation of the dc program (to be introduced subsequently).
In particular, the class of good dc functions
will play an important role in two algorithms that we introduce later; see Propositions~\ref{pr:penalized dc}
and \ref{pr:convergence distributed}.  Incidentally, since every quadratic function is a differentiable dc function,
it follows that every convex constrained optimization problem with a quadratic objective is a ``good dc program'' 
while remaining possibly nonconvex.

\subsection{Concepts of stationarity}

As a non-convex, non-differentiable optimization program, there are many kinds of stationary solutions
for a dc program.  Ideally, we want to be able to identify a stationary solution of the sharpest kind.  Arguably,
for the convex constrained dc program (\ref{eq:ccdc}), a d(irectional)-stationary solution defined
in terms of the directional derivatives of the objective function would qualify for this purpose.  In what follows,
we clarify the relations of several major kinds of stationary solutions of (\ref{eq:general_problem}) by starting
with the definition of d-stationarity.

\gap

Specifically, we say that a vector $x \in X$ is a (constrained) d(irectional)-stationary point of $\zeta$ on $X$
if
\begin{equation} \label{eq:df d-stationary}
\zeta^{\, \prime}(x;x^{\, \prime} - x) \, \geq \, 0, \epc \forall \, x^{\, \prime} \, \in \, X,
\end{equation}
or equivalently, $f^{\, \prime}(x;x^{\, \prime} - x) \geq g^{\, \prime}(x;x^{\, \prime} - x)$ for
all $x \in X$.  Since $g^{\, \prime}(x;d) = \displaystyle{
\max_{v \in \partial g(x)}
} \, v^Td$, where $\partial g(x)$ is the subdifferential of the convex function $g$ at $x$, it follows
that $x$ is a (constrained) d-stationary point of the dc function $\zeta$ on $X$ if for all
$v \in \partial g(x)$,
\begin{equation} \label{eq:dd+subgrad}
f^{\, \prime}(x;x^{\, \prime} - x ) \, \geq \, v^T( \, x^{\, \prime} - x ), \epc \forall \, x^{\, \prime} \, \in X;
\end{equation}
or equivalently, if
\[
x \, \in \, \displaystyle{
\operatornamewithlimits{\mbox{argmin}}_{x^{\, \prime} \in X}
} \, f(x^{\, \prime}) - v^Tx^{\, \prime}, \epc \forall \, v \, \in \, \partial g(x).
\]
Letting $\wh{f} \triangleq f + \delta_X$, where $\delta_X$ is the indicator function of the set $X$, i.e.,
$\delta_X(x) = \left\{ \begin{array}{ll}
0 & \mbox{if $x \in X$} \\
\infty & \mbox{otherwise}
\end{array} \right.$, 
we deduce that $x \in X$ is d-stationary point of $\zeta$ on $X$ if and only if
$v \in \partial \wh{f}(x)$ for all $v \in \partial g(x)$; i.e., if and only if
$\partial g(x) \subseteq \partial \wh{f}(x) = \partial f(x) + {\cal N}(x;X)$, where ${\cal N}(x;X)$ is
the normal cone of the convex set $X$ at $x \in X$ \cite{Rockafellar70}.
This characterization of a d-stationary point is precisely the notion of a {\sl generalized KKT point} employed in the
dc literature \cite{HoiAnTao14-1,HoiAnHuynhTao14,TTAn97} that is convex analysis based.  We prefer
to follow a directional derivative based definition with the constraint set $X$ exposed in the condition (\ref{eq:df d-stationary})
to facilitate the practical solution of the dc program (\ref{eq:ccdc}).  A weaker notion of stationarity, called
{\sl criticality} in the dc literature, is defined by the condition:
$\partial g(x) \, \cap \, \left( \partial f(x) + {\cal N}(x;X) \right) \neq \emptyset$.  In terms of directional derivatives,
this condition says that $x \in X$ is a {\sl critical point} of $\zeta$ of $X$ if {\sl there exists} (as opposed to {\sl for all})
$v \in \partial g(x)$ such that (\ref{eq:dd+subgrad}) holds.

\gap

Using the C-directional derivative, we say that a vector $x \in X$ is {\sl C(larke)-stationary}
if $\zeta^0(x;x^{\, \prime} - x) \geq 0$ for all $x^{\prime} \in X$.  For a good dc function $\zeta$,
d-stationarity and C-stationarity are equivalent.   We will momentarily provide an example to show that if $\zeta$ is not good,
then the converse implication of C-stationarity implying d-stationarity is not always valid.
This example uses the following fact which is by itself of independent interest.

\begin{proposition} \label{pr:good and bad} \rm
Let $\zeta$ be a dc function defined on an open convex set $\Omega \subseteq \mathbb{R}^n$.
The following two statements are equivalent:
\begin{description}
\item[\rm (a)] both $\zeta$ and its negative are good on $\Omega$;
\item[\rm (b)] $\zeta$ is strictly differentiable on $\Omega$ and there exists a strictly differentiable
function $v$ on $\Omega$ such that $\zeta + v$ is convex.
\end{description}
\end{proposition}

\noindent {\bf Proof.}  (a) $\Rightarrow$ (b).  It suffices to show that $\zeta$ is strictly differentiable
on $\Omega$ if (a) holds.  Since both $\zeta$ and $-\zeta$ are C-regular, we have
\[ \begin{array}{lll}
-\displaystyle{
\liminf_{\stackrel{y \to x}{\tau \downarrow 0}}
} \, \displaystyle{
\frac{\zeta(y + \tau d) - \zeta(y)}{\tau}
} & = & \displaystyle{
\limsup_{\stackrel{y \to x}{\tau \downarrow 0}}
} \, \displaystyle{
\frac{-\zeta(y + \tau d) + \zeta(y)}{\tau}
} \, = \, ( -\zeta )^0(x;d) \\ [0.2in]
& = & ( -\zeta )^{\, \prime}(x;d) \, = \, -\displaystyle{
\lim_{\tau \downarrow 0}
} \, \displaystyle{
\frac{\zeta(x + \tau d) - \zeta(x)}{\tau}
} \, = \, - \displaystyle{
\limsup_{\stackrel{y \to x}{\tau \downarrow 0}}
} \, \displaystyle{
\frac{\zeta(y + \tau d) - \zeta(y)}{\tau}.
}
\end{array} \]
Hence,
\[ \displaystyle{
\liminf_{\stackrel{y \to x}{\tau \downarrow 0}}
} \, \displaystyle{
\frac{\zeta(y + \tau d) - \zeta(y)}{\tau}
} \, = \, \zeta ^{\, \prime}(x;d) \, = \, \displaystyle{
\limsup_{\stackrel{y \to x}{\tau \downarrow 0}}
} \, \displaystyle{
\frac{\zeta(y + \tau d) - \zeta(y)}{\tau} .
}
\]
Consequently,
\[
\displaystyle{
\lim_{\stackrel{y \to x}{\tau \downarrow 0}}
} \, \displaystyle{
\frac{\zeta(y + \tau d) - \zeta(y)}{\tau}
} \, = \, \zeta^{\, \prime}(x;d), \epc \forall \, x \, \in \, \Omega \mbox{ and } \forall \, d \, \in \, \mathbb{R}^n.
\]
Using this limit, we show that $\zeta^{\, \prime}(x;\bullet)$ is linear on $\mathbb{R}^n$ for fixed $x$.
Indeed, we have, for any $d$ and $d^{\, \prime}$ in $\mathbb{R}^n$,
\[ \begin{array}{lll}
\zeta^{\, \prime}(x;d+d^{\, \prime}) - \zeta^{\, \prime}(x;d)
& = & \displaystyle{
\lim_{\stackrel{y \to x}{\tau \downarrow 0}}
} \, \left[ \, \displaystyle{
\frac{\zeta(y + \tau d + \tau d^{\, \prime}) - \zeta(y)}{\tau}
} - \displaystyle{
\frac{\zeta(y + \tau d) - \zeta(y)}{\tau}
} \, \right] \\ [0.3in]
& = & \displaystyle{
\lim_{\stackrel{y \to x}{\tau \downarrow 0}}
} \, \displaystyle{
\frac{\zeta(y + \tau d + \tau d^{\, \prime}) - \zeta(y + \tau d)}{\tau}
} \\ [0.3in]
& = & \displaystyle{
\lim_{\stackrel{y^{\, \prime} \to x}{\tau \downarrow 0}}
} \, \displaystyle{
\frac{\zeta(y^{\, \prime} + \tau d^{\, \prime}) - \zeta(y^{\, \prime})}{\tau}
} \, = \, \zeta^{\, \prime}(x;d^{\, \prime}).
\end{array} \]
The strict differentiability of $\zeta$ follows readily from a direct verification of this
property.  

\gap

(b) $\Rightarrow$ (a).  This follows easily from the trivial equality $\zeta = ( \zeta + v ) - v$.  \hfill $\Box$

\gap

The example below shows that for a dc function whose negative is good, a C-stationary
point is not necessarily d-stationary. % thus cannot be a local minimizer of the function.

\gap

{\bf Example}.  Consider the univariate dc function $\zeta(x) \triangleq 1 + x^2 - 2| x |$
in the scalar variable $x$.  Since $-\zeta$ is a good dc function and $\zeta$
is not differentiable at $x = 0$, $\zeta$ cannot be good.  Clearly,
$\partial_C \zeta(0) = [ -2, 2 ]$ contains the origin; thus $x = 0$ is a C-stationary point.
Yet, $\zeta^{\, \prime}(0;\pm 1) = -2$; thus $x = 0$ is not d-stationary.  \hfill $\Box$

\subsection{Lifted stationarity $\Leftrightarrow$ weak d-stationary} \label{subsec:lifted}

A certain class of nonsmooth dc programs can be ``lifted'' to become a smooth, albeit still nonconvex, program
to which standard stationarity conditions can be applied.  Specifically, consider a dc function of the following kind:
\begin{equation} \label{eq:zeta}
\zeta(x) \, \triangleq \, \phi(x) - \displaystyle{
\max_{\mu \in {\cal M}}
} \, \psi(x,\mu),
\end{equation}
where $\phi$ is a convex function, $\psi$ is convex-concave, i.e., $\psi(\bullet,\mu)$ is convex and $\psi(x,\bullet)$ is concave,
and ${\cal M}$ is a compact set in $\mathbb{R}^{\ell}$.  By not requiring ${\cal M}$ to be convex allows us to include the
case where ${\cal M}$ is a discrete set such as $P \cap \{ 0, 1 \}^{\ell}$, where $P$ is a polyhedron in $\mathbb{R}^{\ell}$,
so that the $\mu$-maximization problem corresponds to a binary optimization problem.  In the event that
$\psi(x,\bullet)$ is linear for fixed $x$, the maximization
of $\psi(x,\mu)$ for $\mu$ in a discrete set is equivalent to the maximization
of $\psi(x,\mu)$ for $\mu$ in the convex hull of the set.  If $x$ is a d-stationary point of $\zeta$ on $X$, then
by the renowned Danskin's Theorem,
\[
\phi^{\, \prime}(x;x^{\, \prime} - x) - \displaystyle{
\max_{\mu \in {\cal M}(x)}
} \, \psi(\bullet,\mu)^{\, \prime}(x;x^{\, \prime} - x ) \, \geq \, 0, \epc \forall \, x^{\, \prime} \, \in \, X,
\]
where ${\cal M}(x) \triangleq \displaystyle{
\operatornamewithlimits{\mbox{argmax}}_{\mu \in {\cal M}}
} \, \psi(x,\mu)$.  Equivalently,
\begin{equation} \label{eq:d-stationary minmin}
\phi^{\, \prime}(x;x^{\, \prime} - x) \, \geq \, \psi(\bullet,\mu)^{\, \prime}(x;x^{\, \prime} - x ),
\epc \forall \, x^{\, \prime} \, \in \, X \mbox{ and } \forall \, \mu \, \in \, {\cal M}(x).
\end{equation}
We say that a vector $x \in X$ is a {\sl weak d-stationary point} of $\zeta$ given by (\ref{eq:zeta}) on $X$
if there exists $\mu \in {\cal M}(x)$ such that
\[
\phi^{\, \prime}(x;x^{\, \prime} - x) \, \geq \, \psi(\bullet,\mu)^{\, \prime}(x;x^{\, \prime} - x ),
\epc \forall \, x^{\, \prime} \, \in \, X.
\]
In contrast, since $\partial \, \displaystyle{
\max_{\mu \in {\cal M}}
} \, \psi(x,\mu) = \mbox{ convex hull of $\partial_x \psi(x,\mu)$ for $\mu \in {\cal M}(x)$}$, where
$\partial_x \psi(x,\mu)$ denotes the subdifferential of the function $\psi(\bullet,\mu)$, it follows
that $x$ is a critical point if there exist finitely many $\mu^i \in {\cal M}(x)$ for $i = 1, \cdots I$,
finitely many nonnegative
scalars $( \lambda_i )_{i=1}^I$ summing to unity, and subgradients $v^i \in \partial_x \psi(x,\mu^i)$ such that
\[
\displaystyle{
\sum_{i=1}^I
} \, \lambda_i \, v^i \, \in \, \partial \phi(x) + {\cal N}(x;X),
\]
which is equivalent to
\begin{equation} \label{eq:critical in lifted}
\phi^{\, \prime}(x;x^{\, \prime} -x ) \, = \, \displaystyle{
\max_{u \in \partial \phi(x)}
} \, u^T( x^{\, \prime} - x ) \, \geq \, \left[ \, \displaystyle{
\sum_{i=1}^I
} \, \lambda_i \, v^i \, \right]^T( x^{\, \prime} - x ) \, \geq \, 0, \epc \forall \, x^{\, \prime} \, \in \, X.
\end{equation}
The latter inequality confirms that d-stationarity $\Rightarrow$ weak d-stationary $\Rightarrow$ criticality;
the reason for these one-sided implications is twofold: (i) the multiplicity of the argmax ${\cal M}(x)$,
and (ii) the multiplicity of the set $\partial_x \psi(x,\bar{\mu})$ even if ${\cal M}(x)$ is the singleton
$\{ \bar{\mu} \}$.  When both ${\cal M}(x)$ and $\partial_x \psi(x,\bar{\mu})$ are singletons,
then d-stationarity $\Leftrightarrow$ weak d-stationary $\Leftrightarrow$ criticality.  This happens
when $\zeta$ is a good dc function.

\gap

Corresponding to the minimization problem (\ref{eq:ccdc}) which takes
the form
\begin{equation} \label{eq:dc min explicit}
\displaystyle{
\operatornamewithlimits{\mbox{minimize}}_{x \in X}
} \, \left[ \, \phi(x) - \displaystyle{
\max_{\mu \in {\cal M}}
} \, \psi(x,\mu) \, \right],
\end{equation}
is the lifted reformulation in the pair of variables $(x,\mu)$:
\begin{equation} \label{eq:dc min lifted}
\displaystyle{
\operatornamewithlimits{\mbox{minimize}}_{ ( x,\mu ) \in X \times {\cal M}}
} \, \left[ \, \phi(x) - \psi(x,\mu) \, \right].
\end{equation}
% where the right-hand minimization in $(x,\mu)$ jointly is the {\sl lifted} formulation of the minimization
% of $\zeta$ on $X$.
In the case where both $\phi$ and $\psi$ are differentiable, the latter minimization has the advantage
over the former in that it is a differentiable program in the variables $(x,\mu)$ jointly, whereas with $\mu$
hidden in the function
$\zeta$, the minimization of $\zeta$ over the $x$-variable alone is not a differentiable program unless
${\cal M}(x)$ is a singleton for all $x$ of interest.

\gap

In general, if ${\cal M}$ is also convex, and $\psi$ is directionally differentiable in both variables jointly (e.g.,
$\psi$ is continuously differentiable in $(x,\mu)$) such that
the total directional derivative is the sum of the partial directional derivatives with respect to the two
arguments, i.e., suppose that
\[
\psi^{\, \prime}((x,\mu);(x^{\, \prime} - x,\mu^{\, \prime} - \mu)) \, = \,
\psi(\bullet,\mu)^{\, \prime}(x;x^{\, \prime} - x) + \psi(x,\bullet)^{\, \prime}(\mu;\mu^{\, \prime} - \mu),
\]
then it is not difficult to show that $(x,\bar{\mu})$ is a stationary point of the function
$\phi(x) - \psi(x,\mu)$ on $X \times {\cal M}$ if and only if $\bar{\mu} \in {\cal M}(x)$ and
\begin{equation} \label{eq:lifted stationary}
\phi^{\, \prime}(x;x^{\, \prime} - x) \, \geq \, \psi(\bullet,\bar{\mu})^{\, \prime}(x;x^{\, \prime} - x ),
\epc \forall \, x^{\, \prime} \, \in \, X.
\end{equation}
Thus, $x$ is a weak d-stationary point of $\zeta$ (given by (\ref{eq:zeta})) on $X$ if and only if
there exists $\bar{\mu} \in {\cal M}(x)$ such that $(x,\bar{\mu})$ is a stationary point of the bivariate function
$\phi(x) - \psi(x,\mu)$ on $X \times {\cal M}$.  In this sense, we can say that $x$ is a {\sl lifted stationary point}
of $\zeta$ on $X$ after we have exposed the $\mu$-variable that is part of the bivariate function $\psi(x,\mu)$.

\gap

The lifted problem (\ref{eq:dc min lifted}) in the joint variables $(x,\mu)$ can be interpreted
as a 2-person Nash equilibrium problem.  Indeed, consider two optimization problems:
one is a minimization problem in the $x$-variable parameterized by $\mu$ and the other is a
maximization in the $\mu$-variable parameterized by $x$:
\begin{equation} \label{eq:player x+mu} \displaystyle{
\operatornamewithlimits{\mbox{minimize}}_{x \in X}
} \, \phi(x) - \psi(x,\mu) \epc \mbox{and} \epc
\displaystyle{
\operatornamewithlimits{\mbox{maximize}}_{\mu \in {\cal M}}
} \, \psi(x,\mu).
\end{equation}
A {\sl Nash equilibrium} \label{FPang09,PRazaviyayn14} of (\ref{eq:player x+mu}) is
a pair $(x^*,\mu^*)$ such that
\[
x^* \, \in \, \displaystyle{
\operatornamewithlimits{\mbox{argmin}}_{x \in X}
} \, \phi(x) - \psi(x,\mu^*) \epc \mbox{and} \epc
\mu^* \, \in \, \displaystyle{
\operatornamewithlimits{\mbox{argmax}}_{\mu \in {\cal M}}
} \, \psi(x^*,\mu).
\]
Since, $\phi - \psi(\bullet,\mu)$ is not necessarily a convex function, we say that
$(x^*,\mu^*)$ is a {\sl quasi-Nash equilibrium} (QNE) \cite{PRazaviyayn14,PScutari11} if
$x^*$ is a stationary point of the differentiable program
\[
\displaystyle{
\operatornamewithlimits{\mbox{minimize}}_{x \in X}
} \, \phi(x) - \psi(x,\mu^*)
\]
and $\mu^* \in \displaystyle{
\operatornamewithlimits{\mbox{argmax}}_{\mu \in {\cal M}}
} \, \psi(x^*,\mu)$.  It is then easy to see that if a pair $(x^*,\mu^*)$ is a QNE of the pair of programs
(\ref{eq:player x+mu}), then $x^*$ is a lifted stationary point of $\zeta$ (given by (\ref{eq:zeta})) on $X$.
Conversely, if $x^*$ is such a stationary solution, then
$(x^*,\mu^*)$ is a QNE of the pair of programs (\ref{eq:player x+mu}) for some $\mu^*$.
The upshot of this observation is that a dc program is intimated related to games through
its equivalent lifted program formulation.

\gap

For the dc minimization problem (\ref{eq:dc min explicit}) and its lifted formulation (\ref{eq:dc min lifted})
with ${\cal M}$ convex, we have the following string of implications that relates different concepts of stationarity.
% extends those in the special case of
% the problem (\ref{eq:general_problem}) by including the notions of C-stationarity and l-stationarity,
% where l(ifted)-stationarity is defined with reference to the joint minimization
% problem (\ref{eq:dc min lifted}) with $\psi(x,\mu)$ being continuously differentiable in $(x,\mu)$ and ${\cal M}$
% is a convex set:
\begin{center}
{\fbox{
\parbox[c]{17.3cm}{
\begin{center}
\begin{tabular}{ccccc}
local minimizer of (\ref{eq:dc min explicit}) & $=====>$ & d-stationary & $=====>$
& lifted stationary \\
& & & & $\Downarrow$ \\
C-stat. $\Leftarrow$ d-stationary & $\stackrel{\mbox{${\cal M}(x)$ singleton}}{<=====}$ & weak d-stationary
& $\stackrel{\mbox{$\psi(x;\bullet)$ linear}}{<======}$ & critical \\
& & & \mbox{$\psi(\bullet;\mu)$ diff.} & \\
& & $\Updownarrow$ & & $\Uparrow$ \\
QNE & $<=====>$ & lifted stationary & $\stackrel{\mbox{$\phi$ diff}}{======>}$ & C-stationary \\
& & & \mbox{$\psi(\bullet;\mu)$ diff} &
\end{tabular}
\end{center}
}}}
\end{center}
Two of the above implications are not accounted for in the above discussion;
namely, criticality implies lifted stationarity if $\psi(x;\bullet)$ is linear and
$\psi(\bullet;\mu)$ is differentiable on $\Omega$ for all $\mu \in {\cal M}$,
and lifted stationarity implies C-stationarity if $\phi$ and $\psi(\bullet,\mu)$
are both differentiable.
To prove the former, let $\{ \mu^i, \lambda_i, v^i \}_{i=1}^I$ be as given in the derivation
of (\ref{eq:critical in lifted}).  Since $\psi(x,\bullet)$ is linear, it follows that
$\displaystyle{
\sum_{i=1}^I
} \, \lambda_i v^i \in \partial_x \psi\left(x,\displaystyle{
\sum_{i=1}^I
} \, \lambda_i u^i \right)$.  Thus, $\displaystyle{
\sum_{i=1}^I
} \, \lambda_i v^i = \nabla_x \psi\left(x,\displaystyle{
\sum_{i=1}^I
} \, \lambda_i u^i \right)$.  Since $\bar{\mu} \triangleq \displaystyle{
\sum_{i=1}^I
} \, \lambda_i u^i \in {\cal M}(x)$, weak d-stationarity follows.
To prove the remaining implication, we recall that
C-stationarity of a vector $\wh{x} \in X$ means that
% $0 \in \partial_C \zeta( \wh{x} ) + {\cal N}_X( \wh{x} )$, or equivalently,
$\zeta^0(\wh{x};x - \wh{x} ) \geq 0$ for all $x \in X$.  By the definition of the C-generalized gradient,
we have, for any vector $d$,
\[ \begin{array}{lll}
\zeta^0(\wh{x};d ) & = & \displaystyle{
\limsup_{\stackrel{y \to \wh{x}}{\tau \downarrow 0}}
} \, \displaystyle{
\frac{\phi(y + \tau d) - \phi(y) - ( \, \varphi(y+ \tau d) - \varphi(y) \, )}{\tau}
} \\ [0.2in]
& \geq & \displaystyle{
\limsup_{\tau \downarrow 0}
} \, \displaystyle{
\frac{\phi( \wh{x} ) - \phi( \wh{x} - \tau d) - ( \, \varphi( \wh{x} ) - \varphi( \wh{x} - \tau d) \, )}{\tau} .
}
\end{array} \]
Let $\mu \in {\cal M}( \wh{x} )$ be such that $\phi^{\, \prime}( \wh{x};x^{\, \prime} - \wh{x}) \geq
\psi(\bullet,\mu)^{\, \prime}( \wh{x};x^{\, \prime} - \wh{x} )$ for all $x^{\, \prime} \, \in \, X$.   We then have,
provided that $\phi$ and $\psi(\bullet,\mu)$ are both differentiable at $\wh{x}$,
\[ \begin{array}{lll}
\zeta^0(\wh{x};x - \wh{x} ) & \geq & \displaystyle{
\limsup_{\tau \downarrow 0}
} \, \displaystyle{
\frac{\phi( \wh{x} ) - \phi( \wh{x} - \tau ( x - \wh{x} )) + \psi( \wh{x} - \tau ( x - \wh{x} ),\mu) - \psi( \wh{x},\mu)}{\tau}
} \\ [0.2in]
& \geq & \nabla \phi( \wh{x} )^T( \, x - \wh{x} \, ) - \nabla_x \psi( \wh{x} )^T( \, x - \wh{x} \, ) .
\end{array} \]
If the value function $\varphi(x)$ is strictly differentiable (thus $\zeta = \phi - \varphi$ is a
good dc function), then all the stationarity concepts discussed so far are equivalent.
% \begin{center}
% {\fbox{
% \parbox[c]{17cm}{
% for a good dc program:
% \begin{center}
% \begin{tabular}{ccccccc}
% d-stationarity & $\Leftrightarrow$ & weak d-stationarity & $\Leftrightarrow$ & criticality & $\Leftrightarrow$ & C-stationarity \\
% & & $\Updownarrow$ & & & & \\
% & & l-stationarity & & & & \\
% & & with ${\cal M}$ convex & & & &
% \end{tabular}
% \end{center}
% }}}
% \end{center}
When ${\cal M}$ is a finite set, $\varphi(x)$ is a piecewise smooth function; its strict differentiability
has been characterized in terms of the gradients $\nabla_x \psi(x,\mu)$ at the maximizing $\mu$'s; see
\cite{QTseng07}.  This class of dc programs, which can be good or not, will be the focus
of our subsequent algorithmic development.

\gap

{\bf Counterexamples}.  We make two remarks with regard to the above string of implications:

\gap

(1) In general, a critical
point of (\ref{eq:dc min explicit}) is not necessarily weakly d-stationary; a counterexample
is provided by the univariate function: $\zeta(x) \triangleq -| x |$ obtained by letting
$\phi(x) \triangleq 0$, $\psi(x,\pm 1) \triangleq \pm x$, ${\cal M} \triangleq \{ \pm 1 \}$,
and $X \triangleq [ -1,1 ]$ for simplicity.  Since $\partial_C \zeta(0) = [ -1,1 ]$
and ${\cal N}_X(0) = \{ 0 \}$, it follows
that $0 \in  \partial_C \zeta(0) + {\cal N}_X(0)$.  Yet
$\displaystyle{
\frac{\partial \psi(0,\pm 1)}{\partial x}
} = \pm 1 \neq 0$.

\gap

(2) If $\phi$ is not differentiable, then a
weak d-stationary point is not necessarily C-stationary.  Take $\phi(x) \triangleq x + | x |$
and the same $\psi$, ${\cal M}$, and $X$ as above, resulting in $\varphi(x) = | x |$; thus $\zeta(x) = x$.
Clearly, $0$ is not a C-stationary point.  Yet, with $\mu = 1$, we have
$\phi^{\, \prime}(0;d) - \displaystyle{
\frac{\partial \psi(0,1)}{\partial x}
} d = d + | d | - d = | d | \geq 0$ for all $d \in \mathbb{R}$.  Hence 0 is a weak d-stationary point;
yet this point has no minimizing property whatsoever with regard to the problem of minimizing $\zeta(x)$
on $X$.  \hfill $\Box$

\gap

Derived from the above discussion, particularly from the counterexamples, the following conclusions
refine our understanding of dc programs and add insights to the existing literature of this
class of non-convex optimization problems.

\gap

$\bullet $ The class of good dc programs, i.e., convex constrained programs whose objectives
are good dc functions, is a favorable class of nonsmooth dc problems
for which many advanced concepts of stationarity are equivalent to the basic d-stationarity that is easily
described and understood in terms of the elementary directional derivatives.

\gap

$\bullet $ Given a dc function (even a differentiable one), a ``bad'' representation as the difference
of two convex functions can yield a weak d-stationary point that is not C-stationary.

\gap

$\bullet $ For general nonsmooth minimization problems, the search for a sharp notion of stationarity has always
been a challenge.  Ideally, one wants to be able to design an algorithm that will compute a stationary
point that has the best chance to be a local minimum. For the class of dc minimization problems
exemplified by (\ref{eq:dc min explicit}), the above examples show that the critical points, C-stationary
points, and even weak d-stationary points are not ideal because it is less likely for them to correspond
to local minima.

\section{dc Constrained dc Programs}

In this section, we study the B-stationarity concept (to be defined momentarily) associated with a
general dc program, i.e., a dc program subject to dc constraints:
\begin{equation} \label{eq:dc program multiple dc constraints}
\begin{array}{ll}
\displaystyle{
\operatornamewithlimits{\mbox{minimize}}_{x \in X}
} & \zeta(x) \, \triangleq \, \phi(x) - \varphi(x) \\ [0.1in]
% & \mbox{where } \varphi(x) \, \triangleq \, \displaystyle{
% \operatornamewithlimits{\max}_{1 \leq i \leq \ell}
% } \, \psi_i(x) \\ [0.2in]
\mbox{subject to} & \phi_{{\rm c},j}(x) - \varphi_{{\rm c},j}(x) \, \leq \, 0, \epc
j \, = \, 1, \cdots, J,
% & \mbox{where } \varphi_{{\rm c}}(x) \, \triangleq \, \displaystyle{
% \operatornamewithlimits{\max}_{1 \leq j \leq \ell_{\rm c}}
% } \, \psi_{{\rm c},j}(x),
\end{array}
\end{equation}
where $\phi$, $\varphi$, $\phi_{{\rm c},j}$, and $\varphi_{{\rm c},j}$ for all $j$ are all convex
functions defined on the open convex set $\Omega$ containing the closed convex set $X$.
This study is not only interesting for its own sake but the results are needed subsequently
in the convergence analysis of an iterative scheme for solving the problem.
Before proceeding, we mention a variation of the problem (\ref{eq:SRsingleProblem}) that
leads to dc constraints; see \cite{Alvarado14}.   Namely, in this version of the problem,
we impose some Quality-of-Service (QoS) constraints defined by a prescribed level
of minimum secrecy rate profile $\vect s^* \triangleq \left( s_q^* \right)_{q=1}^Q$
that need to be satisfied in the power allocation.  Specifically, the problem is
\[ % \begin{equation}\label{eq:SRsingleProblemQoS}
\begin{array}{ll}
\displaystyle{
{\operatornamewithlimits{\mbox{maximize}}_{\left( \mathbf{p}, \mathbf{\wh{p}} \right) \, \geq \, \mathbf{0}}}
} & \displaystyle{
\sum_{q=1}^Q
} \, \displaystyle{
\sum_{k=1}^N
} \, \left[ \, R_{qqk}(\mathbf{p},\mathbf{\wh{p}}) - R_{q0k}(\mathbf{p},\mathbf{\wh{p}}) \, \right]^+ \\ [0.2in]
\mbox{subject to:} & \displaystyle{
\sum_{k=1}^{N}
} \, p_q(k) \, \leq \,  P_q^{\max}  \,\,\,\,\,\,\,\, \forall \, q=1,\ldots,Q \,\,\,\,\,\,\,\ \mbox{(private constraints) } \\ [0.2in]
& \displaystyle{
\sum_{k=1}^{N}
} \, \wh{p}_j(k) \, \leq \, \wh{P}_j^{\max} \,\,\,\,\,\,\,\, \forall \,  j=1,\ldots,J \,\,\,\,\,\,\,\ \mbox{ (coupling constraints) }
\\ [0.2in]
\mbox{and} & \displaystyle{
\sum_{k=1}^N
} \, \left[ \, R_{qqk}(\mathbf{p},\mathbf{\wh{p}}) - R_{q0k}(\mathbf{p},\mathbf{\wh{p}}) \, \right]^+ \, \geq \, s_q^*,
\epc q \, = \, 1, \cdots, Q \hspace{0.1in} \mbox{(QoS constraints)}.
\end{array} \]
% \end{equation}
Since each term $\left[ R_{qqk}(\mathbf{p},\mathbf{\wh{p}}) - R_{q0k}(\mathbf{p},\mathbf{\wh{p}}) \right]^+$ is a dc function
of the power variables, the QoS constraints are of the dc type.  Another class of problems that leads to a dc constrained dc program
is the class of quadratic programs with (linear) complementarity constraints (QPCC) \cite{BMPang13,BMPang14,JRalph99}.  Specifically,
consider
\begin{equation} \label{eq:QPCC}
\begin{array}{ll}
\displaystyle{
\operatornamewithlimits{\mbox{minimize}}_{(x,y) \, \in \, Z}
} & q(x,y) \\ [0.1in]
\mbox{subject to} & 0 \, \leq \, y \, \perp \, r + Nx + My \, \geq \, 0,
\end{array} \end{equation}
where $q(x,y)$ is a (possibly nonconvex) quadratic function, $Z$ is a polyhedron in $\mathbb{R}^{n+m}$,
$r$ is an $m$-dimensional vector, $N$ is an $m \times n$ matrix, $M$ is an $m \times m$ matrix (not
necessarily positive semidefinite), and the $\perp$ notation denotes the complementarity between the
variables $y$ and $w \triangleq r + Nx + My$.  Since the (linear) complementarity constraint is clearly
equivalent to 3 conditions, 2 linear and 1 quadratic:
\[
y \, \geq \, 0, \epc r + Nx + My \, \geq \, 0, \epc \mbox{and} \epc y^T( \, r + Nx + My \, ) \, \leq \, 0,
\]
the QPCC is a linearly constrained dc program with one additional dc constraint.

\gap

In general, multiple dc constraints can be combined into a single nondifferentiable dc constraint using the
max-function.  Indeed, the $J$ dc constraints in (\ref{eq:dc program multiple dc constraints}) are equivalent
to the single dc constraint:
\[
\displaystyle{
\max_{1 \leq j \leq J}
} \, \left( \, \phi_{{\rm c},j}(x) - \varphi_{{\rm c},j}(x) \, \right) \, \leq \, 0.
\]
Note that
\[
\displaystyle{
\max_{1 \leq j \leq J}
} \, \left( \, \phi_{{\rm c},j}(x) - \varphi_{{\rm c},j}(x) \, \right) \\ [0.1in]
\epc = \, \underbrace{\displaystyle{
\max_{1 \leq j \leq J}
} \, \left( \, \phi_{{\rm c},j}(x) + \displaystyle{
\sum_{j \neq \ell=1}^J
} \, \varphi_{{\rm c},\ell}(x) \, \right)}_{\mbox{$\phi_{\rm c}(x)$}} - \underbrace{\displaystyle{
\sum_{\ell=1}^J
} \, \varphi_{{\rm c},\ell}(x)}_{\mbox{$\varphi_{\rm c}(x)$}} ,
\]
where $\phi_{\rm c}(x)$ and $\varphi_{\rm c}(x)$ are both convex functions with the latter being
(strictly) differentiable if each $\varphi_{{\rm c},\ell}$ is so.  Thus $\zeta_{\rm c}(x) \triangleq \phi_{\rm c}(x) - \varphi_{\rm c}(x)$
is a good dc function if each $\varphi_{{\rm c},\ell}$ is strictly differentiable.  Thus, we restrict the discussion below
to a singly dc constrained dc program:
\begin{equation} \label{eq:dc program dc constraints}
\begin{array}{ll}
\displaystyle{
\operatornamewithlimits{\mbox{minimize}}_{x \in X}
} & \zeta(x) \, \triangleq \, \phi(x) - \varphi(x) \\ [0.1in]
\mbox{subject to} & \zeta_{\rm c}(x) \, \triangleq \, \phi_{\rm c}(x) - \varphi_{\rm c}(x) \, \leq \, 0.
\end{array}
\end{equation}
Due to the last nonconvex constraint: $\phi_{{\rm c}}(x) \leq \varphi_{{\rm c}}(x)$,
the above problem is considerably more complicated than the convex constrained problem (\ref{eq:ccdc}).
For one thing, constraint qualifications (CQs) are needed to yield a constructive description of the stationarity
condition of the problem (\ref{eq:dc program dc constraints}); this is not a trivial task as the dc constraint
is both nondifferentiable and nonconvex.
% Overcoming the complications caused by this constraint,
% we first introduce a restricted stationarity concept tailored to the dc constrained dc program (\ref{eq:dc program dc constraints})
% that is based on convex-constrained optimization and does not in principle require any CQ.
We focus on the well-known concept of stationarity based on the B(oulingand) tangent cone
of a constraint set at a feasible point.  Applied to (\ref{eq:dc program dc constraints}), this concept,
called B-stationarity \cite{Pang07}, pertains to a feasible vector $x^* \in \wh{X}$ satisfying
\begin{equation} \label{eq:dd tangent}
\zeta^{\, \prime}(x^*;d) \, \geq \, 0, \epc \forall \, d \, \in \, {\cal T}_{\wh{X}}(x^*),
\end{equation}
where $\wh{X} \triangleq \left\{ x \in X \mid \phi_{{\rm c}}(x) \leq \varphi_{{\rm c}}(x) \right\}$
is the (nonconvex) feasible set of (\ref{eq:dc program dc constraints}) and ${\cal T}_{\wh{X}}(x^*)$
is the Bouligand tangent cone of $\wh{X}$ at $x^* \in \wh{X}$, i.e., $d \in {\cal T}_{\wh{X}}(x^*)$
if there exist a sequence of vectors $\{ x^k \} \subset \wh{X}$ converging to $x^*$ and a sequence of positive
scalars $\{ \tau_k \}$ converging to 0 such that $d = \displaystyle{
\lim_{k \to \infty}
} \, \displaystyle{
\frac{x^k - x^*}{\tau_k}
}$.  For a nonconvex set such as $\wh{X}$ that involves the nondifferentiable function $\zeta_{\rm c}$,
it is difficult to derive a constructive
description of ${\cal T}_{\wh{X}}(x^*)$.  Thus the B-stationarity condition (\ref{eq:dd tangent}) is hard
to verify in general and no existing computational scheme can compute a B-stationary point for the problem
(\ref{eq:dc program dc constraints}) according to this definition.  Incidentally, B-stationarity reduces
to d-stationarity without the dc constraint; we use the former terminology to highlight the nonconvexity
and nondifferentiability of the dc constraint.

\subsection{A subclass of dc constraints}

Our goal here is to introduce a constraint qualification for the special case of
% We next introduce a
% convex-constraint based stationary concept for the dc program
(\ref{eq:dc program dc constraints}) where
% that serves as a bridge toward the computation of a B-stationary point.   In essence, this alternative definition
% of stationarity relies on the observation that a sufficient condition for (\ref{eq:dd tangent}) to hold
% is that
% \[
% \zeta^{\, \prime}(x^*;x - x^*) \, \geq \, 0, \epc \forall \, x \, \in \, \wh{X}
% \]
% and this sufficient condition is necessary if $\wh{X}$ is convex.  Thus the idea is to restrict stationarity to
% a convex component of the nonconvex set $\wh{X}$.  This can be accomplished when
\begin{equation} \label{eq:psic}
\varphi_{\rm c}(x) \, \triangleq \, \displaystyle{
\operatornamewithlimits{\max}_{1 \leq k \leq L}
} \, \psi_{{\rm c},k}(x)
\end{equation}
is the pointwise maximum of finitely many differentiable convex functions but there
is no structural assumption on $\phi_{\rm c}(x)$ except its convexity.  Under the stipulation, the feasible
set $\wh{X}$ is the union of finitely many convex sets consisting of the ``smooth'' pieces of $\wh{X}$.
Specifically, we have
\[
\wh{X} \, = \, \displaystyle{
\bigcup_{j=1}^L
} \, \wh{X}^j, \epc \mbox{with} \epc
\wh{X}^j \, \triangleq \, \left\{ \, x \, \in \, X \, \mid \, \phi_{\rm c}(x) \, \leq \, \psi_{{\rm c},j}(x) \, \right\}.
\]
The approach below is reminiscent of the study of stationarity for
the class of mathematical programs with complementarity
constraints \cite{LPRalph96,PFukushima99,Scholtes04,SScholtes00}, and more generally, problems
with piecewise smooth constraints.  In particular, the stationarity theory in \cite{Scholtes04}
is in principle applicable to the above representation of the feasible $\wh{X}$.  Yet, by focusing on each individual set
$\wh{X}^j$, we are able to derive a full characterization of the tangent cone of this set
at a feasible point $\bar{x}$ under a pointwise CQ of the Slater type.  For the discussion below to be meaningful,
we assume that $\bar{x} \in \wh{X}$ is such that $\phi_{\rm c}(\bar{x}) = \varphi_{\rm c}(\bar{x})$.  Indeed,
if $\phi_{\rm c}(\bar{x}) < \varphi_{\rm c}(\bar{x})$, then ${\cal T}_{\wh{X}}(\bar{x}) = {\cal T}_X(\bar{x})$ and
there is no need to analyze ${\cal T}_{\wh{X}}(\bar{x})$ further because we assume that ${\cal T}_X(\bar{x})$ is
well behaved.

\gap

We introduce a convex subset of $\wh{X}^j$ by linearizing
the function $\psi_{{\rm c},j}(x)$ at the given point $\bar{x} \in \wh{X}^j$, obtaining a convex subset of
$\wh{X}^j$:
\[
\wh{Y}^j(\bar{x}) \, \triangleq \, \left\{ \, x \, \in \, X \, \mid \, \phi_{\rm c}(x) \, \leq \,
\psi_{{\rm c},j}(\bar{x}) + \nabla \psi_{{\rm c},j}(\bar{x})^T( \, x - \bar{x} \, ) \, \right\}.
\]
Notice that we cannot linearize the function $\phi_{\rm c}(x)$ because we do not assume that it is differentiable;
moreover, the set $\wh{Y}^j(\bar{x})$ depends on the given vector $\bar{x}$ whereas $\wh{X}^j$ does not.
Clearly, ${\cal T}_{\wh{Y}^j(\bar{x})}(\bar{x}) \subseteq {\cal T}_{\wh{X}^j}(\bar{x})$.  The example
below shows that this inclusion is proper in general.

\begin{example} \rm
Consider the convex univariate functions $\phi_{\rm c}(x) = x^4$ and $\psi_{\rm c}(x) = x^2$ so that the
set $\wh{X} \triangleq \{ x \in \mathbb{R} \mid x^4 - x^2 \leq 0 \} = [ -1,1 ]$ is a simple interval.
Let $\bar{x} = 0$.  It follows that $\wh{Y}(0) = \{ 0 \} = {\cal T}_{\wh{Y}(0)}(0)$; yet ${\cal T}_{\wh{X}}(0) = \mathbb{R}$.
For this example, note that $x = 1/2$ is an ``algebraic Slater'' point of $\wh{X}$, i.e., the inequality $x^4 \leq x^2$ holds
strictly at this point.   \hfill $\Box$
\end{example}

We next introduce a convex cone that is a candidate for the tangent cones ${\cal T}_{\wh{Y}^j(\bar{x})}(\bar{x})$
and ${\cal T}_{\wh{X}^j}(\bar{x})$:
\[
\wh{C}^j(\bar{x}) \, \triangleq \, \left\{ \, d \, \in \, {\cal T}_X(\bar{x}) \, \mid \,
\phi_{\rm c}^{\, \prime}(\bar{x};d) \, \leq \, \nabla \psi_{{\rm c},j}(\bar{x})^Td \, \right\},
\]
for $j \in {\cal M}_{\rm c}(\bar{x}) \triangleq \{ k \mid \varphi_{\rm c}(\bar{x}) = \psi_{{\rm c},k}(\bar{x}) \}$.
The result below shows that if the above cone has an element that satisfies the inequality therein strictly,
then the two tangent cones ${\cal T}_{\wh{Y}^j(\bar{x})}(\bar{x})$ and ${\cal T}_{\wh{X}^j}(\bar{x})$ are both
equal to $\wh{C}^j(\bar{x})$.  This result
is the key for us to show the convergence of the iterative algorithm to be presented later for computing a B-stationary point
of the dc constrained dc program (\ref{eq:dc program dc constraints}).

\begin{proposition} \label{pr:equal tangent cones} \rm
Let $\bar{x} \in \wh{X}$ be such that $\phi_{\rm c}(\bar{x}) = \varphi_{\rm c}(\bar{x})$
and let $j \in {\cal M}_{\rm c}(\bar{x})$.
Suppose an element $\bar{d} \in {\cal T}_X(\bar{x})$ exists such that $\phi_{\rm c}^{\, \prime}( \bar{x};\bar{d} ) <
\nabla \psi_{{\rm c},j}(\bar{x})^T \bar{d}$.  Then
\begin{equation} \label{eq:tangent cones equal}
{\cal T}_{\wh{Y}^j(\bar{x})}(\bar{x}) \, = \, {\cal T}_{\wh{X}^j}(\bar{x}) \, = \, \wh{C}^{\, j}(\bar{x}).
\end{equation}
Thus ${\cal T}_{\wh{X}^j}(\bar{x})$ is a closed convex cone.
\end{proposition}

\noindent {\bf Proof.}  We have the following inclusions:
\[
{\cal T}_{\wh{Y}^j(\bar{x})}(\bar{x}) \, \subseteq \, {\cal T}_{\wh{X}^j}(\bar{x}) \, \subseteq \, \wh{C}^{\, j}(\bar{x}),
\]
where the second inclusion can easily be proved as follows.  Let $d \in {\cal T}_{\wh{X}^j}(\bar{x})$ with unit norm be
given.  Clearly $d \in {\cal T}_X(\bar{x})$.  Let $\{ x^k \} \subset \wh{X}^j(\bar{x}) \setminus \{ \bar{x} \}$
be a sequence converging to $\bar{x}$ such that $d = \displaystyle{
\lim_{k \to \infty}
} \, \displaystyle{
\frac{x^k - \bar{x}}{\| x^k - \bar{x} \|}
}$.  Since $\phi_{\rm c}(x^k) \leq \psi_{{\rm c},j}(x^k)$ for all $k$ and $\phi_{\rm c}(\bar{x}) = \psi_{{\rm c},j}(\bar{x})$,
it follows readily that $\phi_{\rm c}^{\, \prime}(\bar{x};d) \leq \nabla \psi_{{\rm c},j}(\bar{x})^Td$.  It remains to show
that $\wh{C}^j(\bar{x}) \subseteq {\cal T}_{\wh{Y}^j(\bar{x})}(\bar{x})$.  We first show that any
$\bar{d} \in {\cal T}_X(\bar{x})$ satisfying $\phi_{\rm c}^{\, \prime}( \bar{x}; \bar{d} ) <
\nabla \psi_{{\rm c},j}(\bar{x})^T \bar{d}$ must belong to ${\cal T}_{\wh{Y}^j(\bar{x})}(\bar{x})$.
Indeed, for any such $\bar{d}$, let $\{ \bar{x}^k \} \subset X \setminus \{ \bar{x} \}$
be a sequence converging to $\bar{x}$ such that $\bar{d} = \displaystyle{
\lim_{k \to \infty}
} \, \displaystyle{
\frac{\bar{x}^k - \bar{x}}{\| \bar{x}^k - \bar{x} \|}
}$.  Since $\phi_{\rm c}^{\, \prime}( \bar{x};d ) = \displaystyle{
\lim_{k \to \infty}
} \, \displaystyle{
\frac{\phi_{\rm c}(\bar{x}^k) - \phi_{\rm c}(\bar{x})}{\| \bar{x}^k - \bar{x} \|}
}$, it follows that for all $k$ sufficiently large,
$\phi_{\rm c}(\bar{x}^k) < \psi_{{\rm c},j}(\bar{x}) + \nabla \psi_{{\rm c},j}(\bar{x})^T( \bar{x}^k - \bar{x} )$.
Thus $\bar{x}^k \in \wh{Y}^{j}(\bar{x})$ for all $k$ sufficiently large.  Hence $\bar{d} \in {\cal T}_{\wh{Y}^j(\bar{x})}(\bar{x})$.
For any $d \in \wh{C}^j(\bar{x})$, $d + \tau \bar{d}$ remains in ${\cal T}_X(\bar{x})$ and satisfies:
$\phi_{\rm c}^{\, \prime}(\bar{x}; d + \tau \bar{d} ) <
\nabla \psi_{{\rm c},j}(\bar{x})^T ( d + \tau \bar{d})$ for all $\tau > 0$, by
the subadditivity and positive homogeneity of the directional derivative $\phi_{\rm c}^{\, \prime}(\bar{x};\bullet)$.
Therefore, $d + \tau \bar{d} \in {\cal T}_{\wh{Y}^j(\bar{x})}(\bar{x})$
for all $\tau > 0$.  Since the tangent cone is a closed set, it follows that $d \in {\cal T}_{\wh{Y}^j(\bar{x})}(\bar{x})$,
establishing the equalities in (\ref{eq:tangent cones equal}).
The last assertion of the proposition follows readily from the closedness and convexity of $\wh{C}^{\, j}(\bar{x})$. \hfill $\Box$

\gap

The following remarks are worth noting.

\gap

$\bullet $ The existence of a vector $\bar{d} \in {\cal T}_X(\bar{x})$ such that $\phi_{\rm c}^{\, \prime}( \bar{x};\bar{d} ) <
\nabla \psi_{{\rm c},j}(\bar{x})^T \bar{d}$ is equivalent to the existence of a vector $\wh{x} \in X$ such that
$\phi_{\rm c}(\wh{x}) < \psi_{{\rm c},j}(\bar{x}) + \nabla \psi_{{\rm c},j}(\bar{x})^T ( \wh{x} - \bar{x} )$.  Such a vector
satisfies $\phi_{\rm c}(\wh{x}) < \psi_{{\rm c},j}(\wh{x})$.

\gap

$\bullet $ For an index $j \in {\cal M}_{\rm c}(\bar{x})$ for which $\psi_{{\rm c},j}$ is an affine function, we have
$\wh{X}^j = \wh{Y}^j(\bar{x})$ for any $\bar{x} \in \wh{X}^j$.  Hence ${\cal T}_{\wh{Y}^j(\bar{x})}(\bar{x}) = {\cal T}_{\wh{X}^j}(\bar{x})$
always holds; nevertheless for the second equality in (\ref{eq:tangent cones equal}) to hold, we still need the existence of
the vector $\bar{d}$ as in Proposition~\ref{pr:equal tangent cones}.  \hfill $\Box$

\gap

We define a Slater concept for a vector $\bar{x}$ satisfying the assumptions of Proposition~\ref{pr:equal tangent cones}.

\begin{definition} \label{df:pointwise Slater} \rm
The {\sl pointwise Slater} CQ is said to hold for the set $\wh{X}$ at a vector $\bar{x} \in \wh{X}$ satisfying
% if
% either one of the following two conditions holds:
% \begin{description}
% \item[\rm (a)] $\phi_{\rm c}(\bar{x}) < \varphi_{\rm c}(\bar{x})$;
% \item[\rm (b)]
$\phi_{\rm c}(\bar{x}) = \varphi_{\rm c}(\bar{x})$ if for every index
$j \in {\cal M}_{\rm c}(\bar{x})$, % for which $\psi_{{\rm c},j}$ is nonlinear,
there exists $\bar{d}^j \in {\cal T}_X(\bar{x})$
satisfying $\phi_{\rm c}^{\, \prime}( \bar{x};\bar{d}^j ) <
\nabla \psi_{{\rm c},j}(\bar{x})^T \bar{d}^j$.  \hfill $\Box$
% \end{description}
\end{definition}
Since $\phi_{\rm c}^{\, \prime}(\bar{x};\bullet)$ is convex, it follows that the pointwise Slater CQ holds at a vector
$\bar{x} \in \wh{X}$ satisfying $\phi_{\rm c}(\bar{x}) = \varphi_{\rm c}(\bar{x})$ if and only if there exists
a single vector $\bar{d} \in {\cal T}_X(\bar{x})$ such that $\phi_{\rm c}^{\, \prime}( \bar{x};\bar{d} ) <
\nabla \psi_{{\rm c},j}(\bar{x})^T \bar{d}$ for all $j \in {\cal M}_{\rm c}(\bar{x})$.
At such a point $\bar{x}$, we have the following string of implications:
\begin{center}
\begin{tabular}{ccccc}
pointwise Slater at $\bar{x}$ & $\Rightarrow$ & set algebraic Slater & $\Rightarrow$ & set topological Slater; \ i.e., \\ [5pt]
$\phi_{\rm c}(\wh{x}) < \psi_{{\rm c},j}(\bar{x}) + \nabla \psi_{{\rm c},j}(\bar{x})^T ( \wh{x} - \bar{x} )$ & $\Rightarrow$ &
$\phi_{\rm c}(\wh{x}) < \psi_{{\rm c},j}(\wh{x})$ & $\Rightarrow$ & $\wh{x}$ interior pt.\ of $\Phi_{{\rm c},j}$ \\ [5pt]
$\Updownarrow$ & & & & \\ [5pt]
$\phi_{\rm c}^{\, \prime}( \bar{x};\bar{d} ) < \nabla \psi_{{\rm c},j}(\bar{x})^T \bar{d}$, & & & &
\end{tabular}
\end{center}
where $\Phi_{{\rm c},j} \triangleq \{ x \mid \phi_{\rm c}(x) \leq \psi_{{\rm c},j}(x) \}$.
Nevertheless, the reverse of each of the two implications
is in general not true, the main reason being the nonconvexity of the set $\Phi_{{\rm c},j}$.
The corollary below is an immediate consequence of Proposition~\ref{pr:equal tangent cones}.

\begin{corollary} \label{co:tangent dc} \rm
If $\bar{x} \in \wh{X}$ satisfies the pointwise Slater CQ, then
\begin{equation} \label{eq:tangent cone union}
{\cal T}_{\wh{X}}(\bar{x}) \, = \, \displaystyle{
\operatornamewithlimits{\bigcup}_{j \in {\cal M}_{\rm c}(\bar{x})}
} \, \wh{C}^j(\bar{x}).
\end{equation}
Hence, ${\cal T}_{\wh{X}}(\bar{x})$ is the union of finitely many closed convex cone.  \hfill $\Box$
\end{corollary}

The next result identifies another situation in which the equalities in (\ref{eq:tangent cone union}) will hold.
We recall that a function $\theta$ is {\sl piecewise affine} on a domain ${\cal D}$ \cite[Definition~4.1.3]{FPang03}
if it is continuous and there exist finitely many affine functions $\{ \theta_i \}_{i=1}^K$ such that
$\theta(x) \in \{ \theta_i(x) \}_{i=1}^K$ for all $x \in {\cal D}$.

\begin{proposition} \label{pr:polyhedral case} \rm
Let $\bar{x} \in \wh{X}$ be such that $\phi_{\rm c}(\bar{x}) = \varphi_{\rm c}(\bar{x})$.
If $X$ is a polyhedron and the (convex) function $\phi_{\rm c}$ is piecewise affine on $X$,
then (\ref{eq:tangent cone union}) holds.
\end{proposition}

\noindent {\bf Proof.}  It suffices to show the inclusion: $\wh{C}^j(\bar{x}) \subseteq {\cal T}_{\wh{Y}^j(\bar{x})}(\bar{x})$
for all $j \in {\cal M}_{\rm c}(\bar{x})$.  Let $d \in \wh{C}^j(\bar{x})$.  Since $X$ is polyhedral, it follows that
$\bar{x} + \tau d \in X$ for all $\tau > 0$ sufficiently small.  Moreover, for all such $\tau$, we have
\[
\phi_{\rm c}(\bar{x} + \tau d) \, = \, \phi_{\rm c}(\bar{x}) + \tau \, \phi_{\rm c}^{\, \prime}(\bar{x};d)
\]
by the piecewise affine property of $\phi_{\rm c}$; see Exercise~4.8.10 in \cite{FPang03}.  From this equality, we easily deduce
that $d \in {\cal T}_{\wh{Y}^j(\bar{x})}(\bar{x})$.  \hfill $\Box$

\gap

Based on the last two results, we can derive the following necessary and sufficient conditions for a
B-stationary point of the program (\ref{eq:dc program dc constraints}).

\begin{proposition} \label{pr:B-stationary} \rm
Let $\bar{x} \in \wh{X}$.  Provided that either $\bar{x}$ satisfies the pointwise Slater CQ or the
assumptions of Proposition~\ref{pr:polyhedral case} hold, the following statements are equivalent:
\begin{description}
\item[\rm (a)] $\bar{x}$ is a B-stationary point of (\ref{eq:dc program dc constraints});
\item[\rm (b)] for every $j \in {\cal M}_{\rm c}(\bar{x})$,
$\zeta^{\, \prime}(\bar{x};d) \geq 0$ for all $d$ in $\wh{C}^{\, j}(\bar{x})$;
\item[\rm (c)] for every $j \in {\cal M}_{\rm c}(\bar{x})$, $\bar{x}$ is a d-stationary
point of $\zeta(x)$ on the convex subset $\wh{Y}^j(\bar{x})$ of $\wh{X}^j$; i.e.,
$\zeta^{\, \prime}(\bar{x};x - \bar{x}) \geq 0$ for all $x \in \wh{Y}^j(\bar{x})$.
\hfill $\Box$
\end{description}
\end{proposition}

Thus, under the assumptions of Proposition~\ref{pr:B-stationary},
% given $\bar{x} \in \wh{X}$ satisfying the pointwise Slater CQ,
checking if $\bar{x}$
is a B-stationary point of (\ref{eq:dc program dc constraints}) can be determined
by showing that $\bar{x}$ is a d-stationary of $| {\cal M}_{\rm c}(\bar{x}) |$ convex-constrained
dc programs (part (c)).  An algorithm for accomplishing the latter task is presented in Section~\ref{sec:extension I}.
Nevertheless, the identification of $\bar{x}$ requires an extension of this algorithm to deal with the dc
constraint.

\section{Computing d-Stationary Points}

We are now ready to discuss the next main topic of this paper, namely the
computation of a d-stationary point of a convex constrained dc program and its extension to a
B-stationary point when there is a dc constraint present in the problem.  The discussion is divided into
3 parts: the first part (Subsection~\ref{subsec:basic algorithm})
deals with the convex-constrained dc program (\ref{eq:ccdc}), the second part
extends the discussion to a problem with a dc constraint, and the third and last part discusses a parallel
implementation when the problem objective function has a sum structure.

\subsection{The basic algorithm} \label{subsec:basic algorithm}

The DCA is a well-known algorithm for solving a dc program.  In its abstract form \cite{ATao05,HThoai99,TTAn97,Tuy87},
the algorithm works with the subgradients of the two convex functions $\phi(x)$ and $\varphi(x)$, taken to be extended
valued, whose difference is the objective function of the problem; the constraints are all embedded in $\phi$ and $\varphi$.
Subsequential convergence to a critical point is proved, among other properties of the algorithm.  Assuming that
$\varphi(x)$ is differentiable,
the paper \cite{SLanckriet12} revisited the DCA and extended it, called the convex-concave procedure (CCCP),
to a problem with dc constraints defined by good dc functions.
Thus the setting in the latter reference pertains to good dc programs.  To motivate the discussion below, we give an
example of a dc function that is not good and show that the limit point obtained by the DCA is not a d-stationary solution.

\begin{example} \label{ex:DCA not d-stationary} \rm
Consider the univariate, unconstrained minimization of the dc function $\thalf x^2 - \max(-x,0)$ whose unique d-stationary
point is $x = -1$.  Choose a positive $x^0$ as the initial iterate. Without regularization, the DCA computes $x^1$ by minimizing
$\thalf x^2$, yielding $x^1 = 0$.  At this point, if the subgradient of the plus-function is incorrectly picked,
the algorithm could stay at the origin
forever.  A better illustration is to consider a regularization of the DCA wherein at each iteration $\nu$, the algorithm
minimizes the regularized function $\thalf x^2 - ( \partial \max(-x,0) |_{x=x^{\nu}} ) ( x - x^{\nu} ) + \thalf ( x - x^{\nu} )^2$.
It is not hard to see that starting at the same positive $x^0$, the regularized DCA generates a sequence of iterates satisfying
the recursive equation $x^{\nu+1} = \thalf x^\nu$
for $\nu = 0, 1, \cdots$, which converges to the non-d-stationary point $x^\infty = 0$.  For this example, it is easy to modify
the DCA so that the unique d-stationary point can be computed; one such modification is that at each iteration $\nu$, we
consider 2 subproblems: (i) minimizing $\thalf x^2 + ( x - x^{\nu} ) + \thalf ( x - x^{\nu} )^2$, and (ii) minimizing
$\thalf x^2 + \thalf ( x - x^{\nu} )^2$, and choose the next iterate to be the minimizer of these two subproblems
that gives a lower value of the original objective function.
We leave it to the reader to verify that this modified procedure will converge to the d-stationary solution of $-1$.  \hfill $\Box$
\end{example}

Consider the dc program:
\begin{equation} \label{eq:dc program}
\displaystyle{
\operatornamewithlimits{\mbox{minimize}}_{x \in X}
} \ \zeta(x) \, \triangleq \, \phi(x) - \varphi(x),
\epc \varphi(x) \, \triangleq \, \displaystyle{
\operatornamewithlimits{\max_{1 \leq i \leq \ell}}
} \, \psi_i(x)
\end{equation}
where $\phi$ and each $\psi_i$ are
convex functions defined on an open convex set $\Omega$
containing the feasible set $X$, which is a closed convex set in $\mathbb{R}^n$.  Moreover, we assume
that each $\psi_i$ is continuously differentiable (C$^1$) on $\Omega$.  Being a
pointwise maximum of finitely many C$^1$ convex functions, $\varphi$ is a convex piecewise smooth
function with directional derivative at a point $x$ along a direction $d \in \mathbb{R}^n$ given by
\[
\varphi^{\, \prime}(x;d) \, = \, \displaystyle{
\operatornamewithlimits{\max}_{i \in {\cal M}(x)}
} \, \nabla \psi_i(x)^Td,
\]
where ${\cal M}(x) \triangleq \displaystyle{
\operatornamewithlimits{\mbox{argmax}}_{1 \leq i \leq \ell}
} \, \psi_i(x)$.  For a given scalar $\varepsilon > 0$, let
${\cal M}_{\varepsilon}(x) \triangleq \left\{ i \mid \psi_i(x) \geq \varphi(x) - \varepsilon \right\}$
which is a superset of ${\cal M}(x)$.  The following result gives necessary and sufficient condition
for a d-stationary point of (\ref{eq:dc program}) that is useful for its computation.

\begin{proposition} \label{pr:d-stationary dc} \rm
A vector $\bar{x} \in X$ is a d-stationary solution of (\ref{eq:dc program}) if and only if
for every $i \in {\cal M}(\bar{x})$, $\bar{x} \in \displaystyle{
\operatornamewithlimits{\mbox{argmin}}_{x \in X}
} \left[ \phi(x) - \nabla \psi_i(\bar{x})^T ( x - \bar{x} ) \right]$, or equivalently,
$\bar{x} = \displaystyle{
\operatornamewithlimits{\mbox{argmin}}_{x \in X}
} \left[ \phi(x) - \nabla \psi_i(\bar{x})^T ( x - \bar{x} ) + \thalf \| x - \bar{x} \|^2 \right]$.
\end{proposition}

\noindent {\bf Proof.} This follows readily because both functions $\phi(x) - \nabla \psi_i(\bar{x})^T ( x - \bar{x} )$
and
$\phi(x) - \nabla \psi_i(\bar{x})^T ( x - \bar{x} ) + \thalf \| x - \bar{x} \|^2$ are convex in $x$.  \hfill $\Box$

\begin{center}
{\fbox{
\parbox[c]{17cm}
{{\bf Algorithm I.}  Let $\varepsilon > 0$ be given.  For a given $x^{\nu} \in X$ at iteration $\nu$ and for each
index $i \in {\cal M}_{\varepsilon} (x^\nu)$, let
%$i_{\nu} \triangleq \iota(x^{\nu})$.
%For each index $i \in {\cal M}_{\varepsilon}(x^{\nu})$, let
\begin{equation} \label{eq:subproblem}
\wh{x}^{\, \nu,i} \, = \, \displaystyle{
\operatornamewithlimits{\mbox{argmin}}_{x \in X}
} \ \phi(x) - \psi_i(x^\nu) - \nabla \psi_{i}(x^{\nu})^T( x - x^{\nu}) +\thalf \, \| \, x - x^{\nu} \, \|^2.
\end{equation}
Let $\wh{i} \in \displaystyle{
\operatornamewithlimits{\mbox{argmin}}_{i \in {\cal M}_{\varepsilon}(x^\nu)}
} \, \zeta(\wh{x}^{\,\nu,i})+ \thalf \|\wh{x}^{\,\nu,i} - x^\nu\|^2$; set $x^{\nu+1} \triangleq \wh{x}^{\,\nu,\wh{i}}$.
If $x^{\nu+1} = x^{\nu}$, terminate; otherwise replace $\nu$ by $\nu + 1$ and repeat the iteration.
}}}
\end{center}

Before proving the convergence of the above algorithm, we offer a few comments.  First of all, with
$\varepsilon = 0$ and noting that $\nabla \psi_i(x^{\nu}) \in \partial \varphi(x^{\nu})$, the algorithm is
the ``complete primal DCA'' described in Section 3 of the unpublished report \cite{HoiAnTao13}.  Nevertheless,
as shown by the example above, convergence to a d-stationary point of the algorithm with $\varepsilon = 0$ cannot be
proved.  Thus, the introduction of the scalar $\varepsilon$ is essential.  Second, the proximal regularization
is perhaps not needed as one can always strongly convexify the functions $\phi$ and $\psi_i$ without changing the
difference function $\zeta$; indeed we always have
\[
\zeta(x) \, = \, \left[ \, \phi(x) + c(x) \, \right] - \displaystyle{
\operatornamewithlimits{\max_{1 \leq i \leq \ell}}
} \, \left[ \, \psi_i(x) + c(x) \, \right],
\]
for any strongly convex function $c(x)$ such as $\thalf \| x \|^2$.  We adopt the term $\thalf \| x - x^{\nu} \|^2$
as this is a common regularization in many nonlinear programming algorithms.  Third, we have left open the practical
solution of the subproblems (\ref{eq:subproblem}) which may require yet an iterative process.  With the many advances
in convex programming in recent years, this is a safe omission as we are adopting this technology as the workhorse in
the above algorithm.  In this regard, one could incorporate a variable step size in the quadratic term instead
of a unit step size to increase flexibility in the practical implementation of the algorithm.  We have omitted all
these refinements as we want to present the basic version of the algorithm and establish its (subsequential) convergence
to a d-stationary point of the program (\ref{eq:dc program}) which we present in the result below.

\begin{proposition} \label{pr:convergence algorithm} \rm
Suppose that the dc function $\zeta$ is bounded below on the closed convex set $X$.
Starting at any $x^0 \in X$ for which the level set
$L(x^0) \triangleq \{ x \in X \mid \zeta(x) \leq \zeta(x^0) \}$ is bounded,
Algorithm I generates a well-defined bounded sequence $\{ x^{\nu} \}$ such that every accumulation
point, at least one of which must exist, is a d-stationary solution of (\ref{eq:dc program}).
Moreover, if the algorithm does not terminate in a finite number of iterations, any such point cannot
be a local maximizer of $\zeta$ on $X$.
\end{proposition}

\noindent {\bf Proof.}  By the update rule of the algorithm, we have
\[
\begin{array}{ll}
\zeta(x^{\nu}) & = \, \phi(x^{\,\nu}) - \displaystyle{
\operatornamewithlimits{\max}_{1 \leq i \leq \ell}
} \, \{\psi_i(x^{\,\nu})\} \\ [0.2in]
& = \, \phi(x^\nu) - \psi_i(x^\nu),\quad \forall i \in {\cal M}(x^\nu) \\ [0.15in]
& \geq \, \phi(\wh{x}^{\,\nu,i}) - \psi_i(x^\nu) - \nabla \psi_i(x^\nu)^T (\wh{x}^{\,\nu,i} - x^\nu) +
\thalf \, \|\wh{x}^{\,\nu,i} - x^\nu\|^2,\quad \forall i \in {\cal M}(x^\nu) \\ [0.15in]
& \epc \mbox{by the definition of $\wh{x}^{\,\nu,i}$} \\ [0.15in]
& \geq \, \phi(\wh{x}^{\,\nu,i}) - \psi_i(\wh{x}^{\,\nu,i}) + \thalf \, \| \, \wh{x}^{\,\nu,i} - x^\nu \, \|^2,
\quad \forall \, i \, \in \, {\cal M}(x^\nu) \\ [0.15in]
& \epc \mbox{by the convexity of $\psi_i$} \\ [0.15in]
& \geq \, \phi(\wh{x}^{\,\nu,i}) - \displaystyle{
\operatornamewithlimits{\max}_{1 \leq j \leq \ell}
} \, \psi_j(\wh{x}^{\,\nu,i}) + \thalf \, \| \, \wh{x}^{\,\nu,i} - x^\nu \, \|^2,
\quad \forall \, i \, \in \, {\cal M}(x^\nu) \\ [0.2in]
& = \, \zeta(\wh{x}^{\,\nu,i}) + \thalf \, \| \, \wh{x}^{\,\nu,i} - x^\nu \, \|^2,
\quad \forall \, i \, \in \, {\cal M}(x^\nu) \\ [0.15in]
& \geq \, \zeta(x^{\nu+1}) +\thalf \, \| \, x^{\nu+1} - x^{\nu} \, \|^2 \epc \mbox{by the definition of $x^{\nu+1}$}.
\end{array} \]
%% &\geq \, \phi(x^{\,\nu+1}) - \psi_i(x^{\,\nu})-\nabla \psi_{i}(x^{\nu})^T( x^{\, \nu+1} - x^{\nu}) +\thalf \, \| \, x^{\, \nu+1} - x^{\nu} \, \|^2 \nonumber\\
%%%&\geq \phi(x^{\,\nu+1}) - \psi_i(x^{\,\nu+1}) +\thalf \, \| \, x^{\, \nu+1} - x^{\nu} \, \|^2\nonumber\\
%%& \geq \theta(\wh{x}^{\,\nu,i}) +\thalf \, \| \, \wh{x}^{\,\nu,i} - x^{\nu} \, \|^2 \nonumber\\
% where the first inequality is due to the definition of $\wh{x}^{\,\nu,i}$ and the second inequality
% is because of the convexity of the function $\psi_i$.
% Hence with $i = \wh{i}$, it follows that
% \[
% \zeta(x^{\nu}) \, \geq \, \zeta(x^{\,\nu+1}) +\thalf \, \| \, x^{\, \nu+1} - x^{\nu} \, \|^2.
% \]
Hence, the sequence of objective values $\{ \zeta(x^{\nu}) \}$ is non-increasing, and strictly decreasing
if % $\displaystyle{
% \operatornamewithlimits{\max_{i \in {\cal M}_{\varepsilon}(x^{\nu})}}
% } \, \|
$x^{\nu+1} \neq x^{\nu}$ for all $\nu$.  Since $\zeta$ is bounded below on $X$, it follows that
$\displaystyle{
\lim_{\nu \to \infty}
} \zeta(x^{\nu})$ exists and
\begin{equation} \label{eq:2limits}
\displaystyle{
\lim_{\nu \to \infty}
} \, \left[ \, \zeta(x^{\nu}) - \zeta(x^{\nu+1}) \, \right] \, = \, \displaystyle{
\lim_{\nu \to \infty}
}  \, \| \, x^{\nu+1} - x^{\nu} \, \| \, = \, 0.
\end{equation}
Since the sequence $\{ x^{\nu} \}$ is contained in the bounded set $L(x^0)$, it has at least one accumulation point.
Let $\{ x^{\nu} \}_{\nu \in \kappa}$ be a subsequence converging to a limit
$x^{\infty}$, which must necessarily belong to $X$.
% It follows from (\ref{eq:2limits}) and the above string of inequalities that
% \[
% \displaystyle{
% \lim_{\nu (\in \kappa) \to \infty}
% } \, \left[ \, \zeta(\wh{x}^{\,\nu,i}) + \thalf \, \| \, \wh{x}^{\,\nu,i} - x^\nu \, \|^2 \, \right] \, = \,
% \zeta(x^{\infty}).
% \]
By restricting the subsequence on hand, a simple limiting argument shows that
${\cal M}(x^{\nu}) \subseteq {\cal M}(x^{\infty}) \subseteq {\cal M}_{\varepsilon}(x^{\nu})$ for all
$\nu \in \kappa$ sufficiently large.
Therefore, using the update rule of the algorithm, for all $i \in {\cal M}(x^{\infty})$, we have
\[ \begin{array}{lll}
\zeta(x^{\nu+1}) +  \thalf \| \, x^{\nu+1} - x^\nu \|^2 & \leq &
\zeta(\wh{x}^{\,\nu,i})+ \thalf \| \, \wh{x}^{\,\nu,i} - x^\nu \, \|^2 \\ [0.1in]
& \leq & \phi(x) - \left( \, \psi_i(x^\nu) + \nabla \psi_i(x^\nu)^T (x-x^\nu) \, \right)
+ \thalf \, \| \, x - x^{\nu} \, \|^2 \epc \forall \, x \, \in \, X.
\end{array}
\]
Taking the limit $\nu (\in \kappa) \rightarrow \infty$ yields
\[
\zeta(x^\infty) \, \leq \, \phi(x) - \left( \, \psi_i(x^\infty) + \nabla \psi_i(x^\infty)^T (x-x^\infty) \, \right)
+ \thalf \, \| \, x - x^{\infty} \, \|^2, \epc \forall \, x \, \in \, X, \ \forall \, i \, \in \, {\cal M}(x^{\infty}),
\]
or equivalently
\[
\phi(x^\infty) \, \leq \phi(x) - \nabla \psi_i(x^\infty)^T (x-x^\infty) + \thalf \, \| \, x - x^{\infty} \, \|^2,
\epc \forall \, x \, \in \, X, \ \forall \, i \, \in \, {\cal M}(x^{\infty}).
\]
The d-stationarity of $x^{\infty}$ for the minimization problem (\ref{eq:dc program}) follows from Proposition~\ref{pr:d-stationary dc}.
To prove the last statement of the proposition, we note that the sequence $\{ \zeta(x^{\nu}) \}$ must be strictly decreasing
(since the algorithm does not terminate in finite number of iterations); moreover, if $\wh{x}$ is any accumumation point of the sequence
$\{ x^{\nu} \}$, then the sequence $\{ \zeta(x^{\nu}) \}$ converges to $\zeta( \wh{x} )$.  Let
$\{ x^{\nu} \}_{\nu \in \kappa^{\, \prime}}$ be a subsequence converging to $\wh{x}$.  We must have
$\zeta(x^{\nu-1}) > \zeta( x^{\nu} ) \geq \zeta( \wh{x} )$ for all $\nu \in \kappa^{\, \prime}$.  Since
$\{ x^{\nu - 1} \}_{\nu \in \kappa^{\, \prime}}$ also converges to $\wh{x}$, it follows that $\wh{x}$ cannot be a local
maximizer of $\zeta$ on $X$.  \hfill $\Box$

\gap

We make several additional remarks regarding the algorithm and its convergence proof.

\gap

$\bullet $ The choice of $\varepsilon > 0$ is important as Example~\ref{ex:DCA not d-stationary} shows the failure of the
algorithm with $\varepsilon = 0$.

\gap

$\bullet $ A major departure of the algorithm from the DCA is that instead of choosing a subgradient from
$\partial \varphi(x^{\nu})$ at iteration $\nu$, we choose the family of
gradients $\{ \nabla \psi_k(x^{\nu}) \}_{k \in {\cal M}_{\varepsilon}(x^{\nu})}$, which is a finite subset
of $\partial \varphi(x^{\nu})$, at the expense of solving multiple convex subprograms.  This extra effort per
iteration leads to the (subsequential) convergence to a d-stationary point of (\ref{eq:dc program}).

\gap

$\bullet $ A referee asked the question of whether the non-increasing property of the sequence of objective
values $\{ \zeta(x^{\nu}) \}$ can be derived from the well-known properties of the proximal map.  This does
not appear to be the case; however, the proof given above is fairly elementary.

\gap

$\bullet $ If the function $\phi = \phi_{\rm nd} + \phi_{\rm d}$ is the sum of two convex functions
with $\phi_{\rm nd}$ being nondifferentiable and $\phi_{\rm d}$ being differentiable, then
we keep $\phi_{\rm nd}$ as it is but approximate $\phi_{\rm d}$ by its first-order Taylor expression
at $x^{\nu}$.  Specifically, we may define $\wh{x}^{\, \nu,i}$ to be the minimizer of
\[
\displaystyle{
\operatornamewithlimits{\mbox{minimize}}_{x \in X}
} \ \phi_{\rm nd}(x) + \phi_{\rm d}(x^{\nu}) +
\nabla \phi_{\rm d}(x^{\nu})^T( \, x - x^{\nu} \, ) + \thalf \, \| \, x - x^{\nu} \, \|^2 -
\psi_i(x^{\nu}) - \nabla \psi_i(x^{\nu})^T( x - x^{\nu} ),
\]
and the same convergence result can be proved.
% This version of the algorithm will be used subsequently
% when we develop a distributed version of the algorithm for an extended version of the
% problem (\ref{eq:dc program}).

\gap

$\bullet $ At this time, we are not able to extend the algorithm to treat the case where $\varphi(x)$ is the value function
of a continuum family of convex functions, i.e., when $\varphi(x) = \displaystyle{
\max_{y \in Y}
} \, \psi(x,y)$ where $\psi(x,\bullet)$ is a concave function and $Y$ is a compact convex set in $\mathbb{R}^m$ for some
positive integer $m$.  It remains an open challenge to develop a practically implementable and provably convergent
algorithm to compute a d-stationary solution of (\ref{eq:dc program}) in this case.

\gap

Proposition~\ref{pr:convergence algorithm} yields the subsequential convergence of Algorithm~I.  There are various additional conditions
under which sequential convergence can be established.  One such condition is the existence of an {\sl isolated} accumulated point
of the sequence; such a point has the property that it is the unique accumulation point of the sequence within a certain neighborhood
of the point.  We formally state this result in the corollary below; its proof follows immediately from \cite[Proposition~~8.3.10]{FPang03}
and is omitted.

\begin{corollary} \label{co:seq convergence} \rm
Under the assumptions of Proposition~\ref{pr:convergence algorithm}, if one of the accumulation points of the
sequence $\{ x^{\nu} \}$  is isolated, then the sequence converges to it.  \hfill $\Box$
\end{corollary}

A referee pointed out that a recent paper \cite{ABSvaiter13} has established
the convergence of the whole sequence (as opposed to subsequences) produced by various classes of
algorithms to a ``critical point'' for a broad class of nonconvex semi-algebraic problems.  It would be interesting
to investigate whether such a sequential convergence result could be established for Algorithm~I applied to the
dc program (\ref{eq:dc program}) with semi-algebraic functions.

\subsection{A randomized version}

When the set ${\cal M}_{\varepsilon}(x^{\nu})$ contains a large number of
elements, then many subproblems (\ref{eq:subproblem}) have to be solved.  Although each of them is convex and presumably easy to solve,
it would be desirable not to solve too many of them in practical implementation.  Randomization could help in
this regard; i.e., we randomize the choice of an appropriate subproblem to be solved at each iteration.
% However, simply
% randomizing the index $i \in {\cal M}_i(x^{\nu})$ at iteration $\nu$ in Algorithm I does not seem to produce a almost surely
% convergent algorithm.  We need to employ the pointwise maximum of 2 (linearized) functions in defining each subproblem.
% For this purpose, we let, for each $x \in X$, $\iota(x)$ be a prescribed index in ${\cal M}(x)$.
We present this randomized algorithm below and show that it will produce a d-stationary
point of the problem (\ref{eq:dc program}) almost surely.

\begin{center}
{\fbox{
\parbox[c]{17cm}
{{\bf The Randomized Version.}  Let a scalar $p_{\min} \in (0,1)$ be given and let $\varepsilon > 0$
be arbitrary.  For a given $x^{\nu} \in X$ at iteration $\nu$, % let $i_{\nu} \triangleq \iota(x^{\nu})$.
choose an index $i \in {\cal M}_{\varepsilon}(x^{\nu})$ randomly so that
\[
p_i^{\nu} \, \triangleq \mbox{ Pr}\left(\mbox{index } i \mbox{ is chosen } \mid x^1, \cdots, x^\nu \, \right)
\, \geq \, p_{\min} > 0.
\]
Let $x^{\nu+1} = \displaystyle{
\operatornamewithlimits{\mbox{argmin}}_{x \in X}
} \ \phi(x) + \thalf \, \| \, x - x^{\nu} \, \|^2 - \left(
% \wh{\zeta}_i(x;x^{\nu})$,
% where % $\wh{\theta}_i(\bullet,x^{\nu})$ is a strongly convex function for given $x^{\nu}$ that is defined by
% \[ \begin{array}{lll}
% \wh{\zeta}_i(x;x^{\nu}) & \triangleq &  \\ [7pt]
% & & \max\left\{ \, \psi_{i_{\nu}}(x^{\nu}) + \nabla \psi_{i_{\nu}}(x^{\nu})^T( x - x^{\nu} ), \,
\psi_i(x^{\nu}) + \nabla \psi_i(x^{\nu})^T( x - x^{\nu} ) \right)$.
% Set
% \[
% x^{\nu+1} \, \triangleq \, \left\{ \begin{array}{ll}
% x^{\nu} & \mbox{if $\zeta(\wh{x}^{\, \nu,i}) + \thalf \| \wh{x}^{\, \nu,i} - x^{\nu} \|^2 > \zeta(x^{\nu})$} \\ [0.1in]
% \wh{x}^{\, \nu,i} & \mbox{otherwise}.
% \end{array} \right.
% \]
}}}
\end{center}

In what follows, we establish the almost sure convergence of the above
randomized algorithm.  For each index $j \in {\cal M}_{\varepsilon}(x^{\nu})$, let
$\wh{x}^{\, \nu,j} \triangleq \displaystyle{
\operatornamewithlimits{\mbox{argmin}}_{x \in X}
} \ \wh{\zeta}_j(x;x^{\nu}) \triangleq \phi(x) + \thalf \, \| x - x^{\nu} \|^2
- \left( \psi_j(x^{\nu}) + \nabla \psi_j(x^{\nu})^T ( x - x^{\nu} ) \right)$.
We have,
\[
\zeta(x^\nu) \, = \, \wh{\zeta}_j(x^\nu;x^\nu)
\, \geq \, \wh{\zeta}_j(\wh{x}^{\, \nu,j};x^\nu) \, \geq \, \zeta(\wh{x}^{\, \nu,j}) +
\thalf \, \| \, \wh{x}^{\, \nu,j} - x^\nu \, \|^2.
\]
Moreover, $x^{\nu+1} = \wh{x}^{\, \nu,j}$ with probability $p_j^{\nu}$.
Taking conditional expectations, the above inequality implies
\[ % \begin{array}{lll}
\mathbb{E}\left[\zeta(x^{\nu+1})\mid x^\nu\right] \, = \, \displaystyle{
\sum_{i \in {\cal M}_{\varepsilon}(x^{\nu})}
} \, p_i^\nu \, \zeta(\wh{x}^{\,\nu,i}) % \\ [0.25in]
\, \leq \, \zeta(x^\nu) - \thalf \, \displaystyle{
\sum_{i \in {\cal M}_{\varepsilon}(x^{\nu})}
} \, p_i^\nu \, \| \, \wh{x}^{\,\nu,i} - x^\nu \, \|^2 .
% \, \leq \, \zeta(x^\nu) - \displaystyle{
% \frac{p_{\min}}{2}
% } \, \displaystyle{
% \sum_{i\in {\cal M}_{\varepsilon}(x^{\nu})}
% } \, \| \, \wh{x}^{\,\nu,i} - x^\nu \, \|^2.
% \end{array}
\]
Consequently, the random sequence $\{\zeta(x^\nu)\}$ is a supermartingale
and assuming that $\zeta$ is bounded from below on $X$, we may conclude that $\{\zeta(x^\nu)\}$ converges
almost surely and, letting $p_i^{\nu} = 0$ for all $i \not\in {\cal M}_{\varepsilon}(x^{\nu})$,
\begin{equation} \label{eq:fixedlimit}
\displaystyle{
\lim_{\nu\rightarrow \infty}
} \, p_i^{\nu} \, \| \, \wh{x}^{\,\nu,i} - x^\nu \, \| \, = \, 0, \;\;\forall \, i \, = \, 1, \cdots, \ell
\end{equation}
with probability one.  In the rest of the proof, we restrict ourselves to the set of probability one in which the above limit
holds.  Consider a point $x^\infty$ that is the limit of the subsequence $\{ x^{\nu} \}_{\nu \in \kappa}$.
By further restricting the subsequence on hand, we can assume that $i_{\nu} = \bar{i}$ for all $\nu \in \kappa$
with $\bar{i}\in \mathcal{M}(x^{\infty})$.  Let $i \in {\cal M}(x^{\infty})$ be given.  It then follows that
$i \in {\cal M}_{\varepsilon}(x^{\nu})$ for all $\nu \in \kappa$ sufficiently large.  Since $p_i^{\nu} \geq p_{\min}$,
it follows from (\ref{eq:fixedlimit}) that
$\displaystyle{
\lim_{\nu (\in \kappa) \rightarrow \infty}
} \, \wh{x}^{\,\nu,i} = \, \displaystyle{
\lim_{\nu (\in \kappa) \rightarrow \infty}
} \, x^\nu \, = \, x^{\infty}$.  Therefore, by the
definition of $\wh{x}^{\,\nu,i}$, we have, for every $x \in X$,
\[
\phi(x) + \thalf \, \| \, x - x^{\nu} \, \|^2 - \psi_i(x^{\nu}) - \nabla \psi_i(x^{\nu})^T( x - x^{\nu} )
\, \geq \, \phi(\wh{x}^{\, \nu,i}) - \psi_i(x^{\nu}) - \nabla \psi_i(x^{\nu})^T( \wh{x}^{\, \nu,i} - x^{\nu} ).
\]
Letting $\nu (\in \kappa) \to \infty$ in the above inequality, we deduce
\[
\phi(x) + \thalf \, \| \, x - x^{\infty} \, \|^2 - \nabla \psi_i(x^{\infty})^T( x - x^{\infty} )\, \geq \, \phi(x^{\infty}),
\]
from which we can deduce that $\phi(x) - \nabla \psi_i(x^{\infty})^T( x - x^{\infty} ) \geq \phi(x^{\infty})$ for all $x \in X$.
By Proposition~\ref{pr:d-stationary dc}, it follows that $x^{\infty}$ is a d-stationary solution of
(\ref{eq:dc program}) almost surely.  This completes the proof of the following convergence result.

\begin{proposition} \label{pr:randomized} \rm
Suppose that the dc function $\zeta$ is bounded below on the closed convex set $X$.
Every limit point of the iterates generated by the randomized algorithm is a d-stationary point of the dc program~\eqref{eq:dc program}
with probability one.  \hfill $\Box$
\end{proposition}

\section{Algorithmic Extension: I} \label{sec:extension I}

In this and the next section, we present two extensions of the deterministic Algorithm I and omit their randomized versions.
When providing the convergence of the extended algorithms, we focus on their subsequential convergence and rely on
Corollary~\ref{co:seq convergence} and the recent reference \cite{ABSvaiter13} for the issue of sequential convergence.
The first extension of Algorithm~I is to the
dc constrained dc program (\ref{eq:dc program dc constraints}).  We start by presenting an immediate consequence of
Propositions~\ref{pr:B-stationary} and \ref{pr:d-stationary dc}.

\begin{proposition} \label{pr:dc constrained dc} \rm
Let $\varphi(x) \triangleq \displaystyle{
\operatornamewithlimits{\max}_{1 \leq j \leq \ell}
} \, \psi_i(x)$ and $\varphi_{\rm c}(x) \triangleq \displaystyle{
\operatornamewithlimits{\max}_{1 \leq j \leq L}
} \, \psi_{{\rm c},j}(x)$, where $\psi_i$ and $\psi_{{\rm c},j}$ are convex differentiable functions on $\Omega$.
Let $\bar{x} \in \wh{X}$ satisfy the pointwise Slater CQ.  It holds that $\bar{x}$ is a B-stationary solution
of (\ref{eq:dc program dc constraints}) if and only if
for every $i \in {\cal M}(\bar{x})$ and every $j \in {\cal M}_{\rm c}(\bar{x})$,
$\bar{x} \in \displaystyle{
\operatornamewithlimits{\mbox{argmin}}_{x \in \wh{Y}^j(\bar{x})}
} \, \left[ \phi(x) - \nabla \psi_i(\bar{x})^T ( x - \bar{x} ) \right]$, or equivalently,
$\bar{x} = \displaystyle{
\operatornamewithlimits{\mbox{argmin}}_{x \in \wh{Y}^j(\bar{x})}
} \, \left[ \phi(x) - \nabla \psi_i(\bar{x})^T ( x - \bar{x} ) + \thalf \| x - \bar{x} \|^2 \right]$.  \hfill $\Box$
\end{proposition}

In the rest of this section, we assume that the two functions $\varphi$ and $\varphi_{\rm c}$ are as
given in Proposition~\ref{pr:dc constrained dc}.
% and that there exists a scalar $\alpha > 0$ such that the set
% \[
% X_{\alpha} \, \triangleq \, \{ \, x \, \in \, X \, \mid \phi_{\rm c}(x) \, \leq \, \varphi_{\rm c}(x) + \alpha \}
% \]
% is nonempty and bounded.
Till now, the issue of feasibility of the problem
(\ref{eq:dc program dc constraints}) has not been addressed.  Indeed, this is a very difficult issue
and we will not directly deal with it.  In what follows, we propose two approaches to
compute a B-stationary point of (\ref{eq:dc program dc constraints}).  The first approach assumes that
a feasible solution of the problem is available which we will use to initiate the algorithm.  The
second approach does not assume that such a (feasible) solution is readily available, perhaps because the
problem is actually not feasible.  We propose a double-loop scheme in which the outer loop solves
a sequence of convex-constrained subproblems by penalizing the dc constraint and the inner loop applies
the basic Algorithm I (or its randomized version) to compute a d-stationary point of the penalized
subproblems.  Convergence of both algorithms will be analyzed.

\subsection{Feasibility assumed}

In this subsection, we assume that a vector $x^0 \in \wh{X}$ is available.  Similar
to the index set ${\cal M}_{\varepsilon}(\bar{x})$ pertaining to the max-function $\varphi(x)$
in the objective,
we define, for each $\varepsilon > 0$ and each $\bar{x} \in X$ the set
% two-step procedure at each iteration.  Given $x^{\nu} \in X$,
% and $t_{\nu} \in [ 0,\alpha ]$ such that $x^{\nu} \in X_{t_{\nu}}$ ,
% we first determine the minimum residual of the constraints $\phi_{\rm c}(x) \leq \psi_{{\rm c},j}(x^{\nu}) +
% \nabla \psi_{{\rm c},j}(x^{\nu})^T (x - x^{\nu})$ for every $j$ in the set
\[
{\cal M}_{{\rm c},\varepsilon}(\bar{x}) \, \triangleq \, \left\{ \, k \, \mid \,
\psi_{{\rm c},k}(\bar{x}) \, \geq \, \varphi_{\rm c}(\bar{x}) - \varepsilon \, \right\}
\]
pertaining to the max-function $\varphi_{\rm c}(x)$ in the constraint.  We also recall the set
\[
\wh{Y}^j(\bar{x}) \, \triangleq \, \left\{ \, x \, \in \, X \, \mid \, \phi_{\rm c}(x) \, \leq \,
\psi_{{\rm c},j}(\bar{x}) + \nabla \psi_{{\rm c},j}(\bar{x})^T( \, x - \bar{x} \, ) \, \right\},
\]
which we have previously defined for a vector $\bar{x} \in \wh{X}^j$ is now extended to
an arbitrary vector $\bar{x} \in X$.  Note: if $\bar{x} \not\in \wh{X}^j$, the nonemptiness of
$\wh{Y}^j(\bar{x})$ is not guaranteed.  Nevertheless $\wh{Y}^j(\bar{x})$ must be nonempty if
$\bar{x} \in \wh{X}$ and $j \in {\cal M}_{\rm c}(\bar{x})$.

\begin{center}
{\fbox{
\parbox[c]{17cm}
{{\bf Algorithm II.}  Let $\varepsilon > 0$ and $x^0 \in \wh{X}$ be given.
At iteration $\nu$, given $x^{\nu} \in \wh{X}$,
we let, for every pair of indices $i \in {\cal M}_{\varepsilon}(x^{\nu})$ and $j \in {\cal M}_{{\rm c},\varepsilon}(x^{\nu})$,
$\wh{x}^{\, \nu,i,j}$ be the (unique) optimal solution of the strongly convex program:
\begin{equation} \label{eq:dc constrained subproblem}
% \begin{array}{ll}
\displaystyle{
\operatornamewithlimits{\mbox{argmin}}_{x \in \wh{Y}^j(x^{\nu})}
} \ \phi(x) - \psi_i(x^{\nu}) - \nabla \psi_i(x^{\nu})^T( x - x^{\nu}) + \thalf \, \| \, x - x^{\nu} \, \|^2
% \\ [0.15in]
% \mbox{subject to} & \phi_{\rm c}(x) \, \leq \,
% \psi_{{\rm c},j}(x^{\nu}) + \nabla \psi_{{\rm c},j}(x^{\nu})^T ( \, x - x^{\nu} \, )
% \end{array}
\end{equation}
if $\wh{Y}^j(x^{\nu}) \neq \emptyset$; otherwise we let $\wh{x}^{\, \nu,i,j} = x^{\nu}$.
Let $\left( \, \wh{i},\wh{j} \, \right) \in \displaystyle{
\operatornamewithlimits{\mbox{argmin}}_{(i,j) \in {\cal M}_{\varepsilon}(x^{\nu}) \times {\cal M}_{{\rm c},\varepsilon}(x^{\nu})}
} \, \zeta(\wh{x}^{\,\nu,i,j})+ \thalf \|\wh{x}^{\,\nu,i,j} - x^{\nu}\|^2$; set $x^{\nu+1} \triangleq \wh{x}^{\,\nu,\wh{i},\wh{j}}$.
}}}
\end{center}

We have the following (subsequential) convergence result of the above algorithm.

\begin{proposition} \label{pr:dc constrained convergence algorithm} \rm
Suppose that the dc function $\zeta$ is bounded below on the feasible set $\wh{X}$.
Starting at any $x^0 \in \wh{X}$ for which the level set
$\wh{L}(x^0) \triangleq \{ x \in \wh{X} \mid \zeta(x) \leq \zeta(x^0) \}$ is bounded,
Algorithm~II generates a well-defined bounded sequence $\{ x^{\nu} \} \subset \wh{X}$
such that every accumulation point $x^{\infty}$, at least one of which must exist,
is feasible to (\ref{eq:dc program dc constraints}); moreover, if $x^{\infty}$ satisfies
the pointwise Slater CQ, then $x^{\infty}$ is a B-stationary point of (\ref{eq:dc program dc constraints}).
\end{proposition}

\noindent {\bf Proof.}  Since $x^{\nu} \in \wh{X}$, it follows that the subproblem (\ref{eq:dc constrained subproblem})
is feasible for all $j \in {\cal M}_{\rm c}(x^{\nu})$ and thus has a unique optimal solution.
Moreover, $x^{\nu+1} \in \wh{X}$ by the gradient
inequality applied to the function $\psi_{{\rm c},\wh{j}}$.
We now follow the proof of Proposition~\ref{pr:convergence algorithm} to deduce the following string
of (in)equalities:
\[
\begin{array}{ll}
\zeta(x^{\nu}) & = \, \phi(x^{\,\nu}) - \displaystyle{
\operatornamewithlimits{\max}_{1 \leq i \leq \ell}
} \, \{\psi_i(x^{\,\nu})\} \\ [0.15in]
& = \, \phi(x^\nu) - \psi_i(x^\nu),\quad \forall i \in {\cal M}(x^\nu) \\ [0.1in]
& \geq \, \phi(\wh{x}^{\,\nu,i,j}) - \psi_i(x^\nu) - \nabla \psi_i(x^\nu)^T (\wh{x}^{\,\nu,i,j} - x^\nu) +
\thalf \, \|\wh{x}^{\,\nu,i,j} - x^\nu\|^2,  \\ [0.15in]
& \epc \forall \, ( i,j ) \in {\cal M}(x^\nu) \times {\cal M}_{\rm c}(x^{\nu});
\mbox{ by the definition of $\wh{x}^{\,\nu,i,j}$} \\ [0.15in]
& \geq \, \phi(\wh{x}^{\,\nu,i,j}) - \psi_i(\wh{x}^{\,\nu,i,j}) + \thalf \, \| \, \wh{x}^{\,\nu,i,j} - x^\nu \, \|^2,
\quad \forall \, ( i,j ) \in {\cal M}(x^\nu) \times {\cal M}_{\rm c}(x^{\nu}) \\ [0.1in]
& \epc \mbox{by the convexity of $\psi_i$} \\ [0.15in]
& \geq \, \phi(\wh{x}^{\,\nu,i}) - \displaystyle{
\operatornamewithlimits{\max}_{1 \leq k \leq \ell}
} \, \psi_k(\wh{x}^{\,\nu,i,j}) + \thalf \, \| \, \wh{x}^{\,\nu,i,j} - x^\nu \, \|^2,
\quad \forall \, ( i,j ) \in {\cal M}(x^\nu) \times {\cal M}_{\rm c}(x^{\nu}) \\ [0.15in]
& = \, \zeta(\wh{x}^{\,\nu,i,j}) + \thalf \, \| \, \wh{x}^{\,\nu,i,j} - x^\nu \, \|^2,
\quad \forall \, ( i,j ) \in {\cal M}(x^\nu) \times {\cal M}_{\rm c}(x^{\nu}) \\ [0.15in]
& \geq \, \zeta(x^{\nu+1}) +\thalf \, \| \, x^{\nu+1} - x^{\nu} \, \|^2 \epc \mbox{by the definition of $x^{\nu+1}$}.
\end{array} \]
As before, it follows that (\ref{eq:2limits}) holds.  Let $\{ x^{\nu} \}_{\nu \in \kappa}$ be a subsequence converging to a limit
$x^{\infty}$, which can easily be seen to belong to $\wh{X}$.  Suppose that $x^{\infty}$ satisfies the pointwise Slater CQ.
According to Proposition~\ref{pr:dc constrained dc}, it suffices to show that $\bar{x} = \displaystyle{
\operatornamewithlimits{\mbox{argmin}}_{x \in \wh{Y}^j(\bar{x})}
} \, \left[ \phi(x) - \nabla \psi_i(\bar{x})^T ( x - \bar{x} ) + \thalf \| x - \bar{x} \|^2 \right]$ for all pairs of indices
$(i,j) \in {\cal M}(x^{\infty}) \times {\cal M}_{\rm c}(x^{\infty})$.  Let $j$ be such an index and
$x \in \wh{Y}^j(x^{\infty})$ be arbitrary.  Then $j \in {\cal M}_{{\rm c},\varepsilon}(x^{\nu})$ for all
$\nu \in \kappa$ sufficiently large.  Let $\bar{x}^j$ satisfy: $\phi_{\rm c}(\bar{x}^j) <
\psi_{{\rm c},j}(x^{\infty}) + \nabla \psi_{{\rm c},j}(x^{\infty})^T( \bar{x}^j - x^{\infty} )$.
For all scalars $\tau \in [ 0,1 )$, with
$x^{\tau} \triangleq \bar{x}^j + \tau ( x - \bar{x}^j ) \in X$, we have
\[
\phi_{\rm c}(x^\tau) \, < \, \psi_{{\rm c},j}(x^{\infty}) + \nabla \psi_{{\rm c},j}(x^{\infty})^T( x^{\tau} - x^{\infty} ).
\]
For each fixed $\tau \in [ 0,1 )$, it follows that for all $\nu (\in \kappa)$ sufficiently large, we have
\[
\phi_{\rm c}(x^{\tau}) \, < \,
\psi_{{\rm c},j}(x^{\nu}) + \nabla \psi_{{\rm c},j}(x^{\nu})^T( \, x^{\tau} - x^{\nu} \, )
\]
for all $\nu \in \kappa$ sufficiently large.
Thus, $x^{\tau}$ is feasible to (\ref{eq:dc constrained subproblem}).  Similar to the proof of
Proposition~\ref{pr:convergence algorithm}, we have for all $i \in {\cal M}(x^{\infty})$,
which is a subset of ${\cal M}_{\varepsilon}(x^{\nu})$ for all $\nu$ sufficiently large
\[ \begin{array}{lll}
\zeta(x^{\nu+1}) +  \thalf \| \, x^{\nu+1} - x^\nu \|^2 & \leq &
\zeta(\wh{x}^{\,\nu,i,j})+ \thalf \| \, \wh{x}^{\,\nu,i,j} - x^\nu \, \|^2 \\ [0.1in]
& \leq & \phi(x^{\tau}) - \left( \, \psi_i(x^\nu) + \nabla \psi_i(x^\nu)^T (x^{\tau}-x^\nu) \, \right)
+ \thalf \, \| \, x^{\tau} - x^{\nu} \, \|^2.
\end{array}
\]
Passing to the limit $\nu (\in \kappa) \to \infty$, we deduce, for all $\tau \in [ 0,1 )$,
\[
\zeta(x^{\infty}) \, \leq \, \phi(x^{\tau}) - \psi_i(x^{\infty}) - \nabla \psi_i(x^\infty)^T ( \, x^{\tau} - x^{\infty} \, )
+ \thalf \| \, x^{\tau} - x^{\infty} \, \|^2.
\]
Passing to the limit $\tau \uparrow 1$, we deduce
\[
\phi(x^{\infty}) \, \leq \, \phi(x) - \nabla \psi_i(x^{\infty})^T( \, x - x^{\infty}) + \thalf \| \, x - x^{\infty} \, \|^2
\]
for all $x \in \wh{Y}^j(x^{\infty})$ and all $(i,j) \in {\cal M}(x^{\infty}) \times {\cal M}_{\rm c}(x^{\infty})$ as desired.
\hfill $\Box$

\subsection{Feasibility not assumed}

Without assuming the feasibility of the problem (\ref{eq:dc program dc constraints}), we propose a penalization of the
dc constraint and establish a limiting result when the penalization tends to infinity.  Penalization techniques in dc programming
and DCA have been investigated for solving
dc constrained dc programs \cite{TTAn97,HoiAnTaoHuynhi12,HoiAnTaoMu99,HoiAnTao14-1}.  The departure of our discussion
from these references is that we aim to compute a B-stationary solution of such a program.  This is accomplished by
considering the following
% Unlike the theory of exact penalization as discussed in \cite[Section~8]{TTAn97},
% which pertains to global minimizers of the penalized subproblems and does not seem appropriate because for all practical
% purposes such minimizers can not be computed, we
% investigate the following stationarity issue that is more relevant due to the nonconvexity.  For a given
% penalty parameter $\rho > 0$, consider the
penalized convex-constrained dc program: for $\rho > 0$,
\[
\displaystyle{
\operatornamewithlimits{\mbox{minimize}}_{x \in X}
} \ \zeta(x) + \rho \, \max\left( \, 0, \zeta_{\rm c}(x) \, \right)
\]
and letting $x^{\rho}$ be a d-stationary point of this problem.  Suppose that for a sequence of penalty parameters $\{ \rho_{\nu} \} \uparrow \infty$,
the corresponding sequence of d-stationary solutions $\left\{ x^{\rho_{\nu}} \right\}$ converges to a limit $x^{\infty}$.  What can we say
about $x^{\infty}$ with regard to the stationarity properties of (\ref{eq:dc program dc constraints})?  Incidentally, since
\[
\phi(x) - \varphi(x) + \rho \, \max\left( \, 0, \phi_{\rm c}(x) - \varphi_{\rm c}(x) \, \right)
\, = \, \left[ \, \underbrace{\phi(x) + \rho \, \max\left( \, \phi_{\rm c}(x), \, \varphi_{\rm c}(x) \, \right)}_{\mbox{convex}} \, \right]
- \left( \, \underbrace{\varphi(x) + \rho \, \varphi_{\rm c}(x)}_{\mbox{convex}} \, \right),
\]
the computation of each $x^{\rho}$ can be accomplished by Algorithm I or its randomized version.  It should be noted that during the computation
of the sequence $\{ x^{\nu} \}$, where $x^{\nu} \triangleq x^{\rho_{\nu}}$, if at any time, an iterate satisfies the dc constraint and is thus
a feasible solution of the problem (\ref{eq:dc program dc constraints}), we have the option of abandoning this penalization approach and return
to the previous direct approach wherein feasibility is maintained throughout the algorithm.  Here, we do not concern ourselves with these algorithmic
details and focus on an understanding of the asymptotic property of the penalization approach employing an unbounded sequence of penalty parameters.
For practical implementation, one should introduce a penalty update rule that circumvents the unboundedness of such a sequence.  Details like this should
best be left for future studies.

\begin{proposition} \label{pr:penalized dc} \rm
In the above setting, the following three statements hold:
\begin{description}
\item[\rm (a)] If $\zeta_{\rm c}(x^{\nu}) > 0$ for infinitely many $\nu$'s, and if $\varphi_{\rm c}$ is
strictly differentiable at $x^{\infty}$, then provided that $X$ is bounded and $\zeta$ is globally Lipschitz continuous on $X$,
$x^{\infty}$ is a d-stationary solution of
\[ \displaystyle{
\operatornamewithlimits{\mbox{minimize}}_{x \in X}
} \ \zeta_{\rm c}(x) \, \triangleq \, \phi_{\rm c}(x) - \varphi_{\rm c}(x);
\]
% i.e., there exists a subgradient $g \in \partial \varphi_{\rm c}(x^{\infty})$ such that
% \[
% \phi_{\rm c}^{\prime}(x^{\infty};x - x^{\infty}) - g^T( \, x - x^{\infty} \, ) \, \geq \, 0, \epc \forall \, x \, \in \, X,
% \]
\item[\rm (b)] If $\zeta_{\rm c}(x^{\nu}) < 0$ for infinitely many $\nu$'s, then $x^{\infty}$ is a d-stationary point of $\zeta$ on $X$,
thus a B-stationary point on $\wh{X}$;
\item[\rm (c)] Suppose $\zeta_{\rm c}(x^{\nu}) = 0$ for all but finitely many $\nu$'s.
If $x^{\infty}$ satisfies the pointwise Slater CQ and ${\cal M}_{\rm c}(x^{\infty})$ is a singleton, then
$x^{\infty}$ is a B-stationary point of (\ref{eq:dc program dc constraints}).
\end{description}
\end{proposition}

{\bf Proof.} (a) If $\zeta_{\rm c}(x^{\nu}) > 0$, then the stationarity condition of $x^{\nu}$ is
\begin{equation} \label{eq:stationary nu}
\zeta^{\prime}(x^{\nu};x - x^{\nu}) + \rho_{\nu} \, \zeta_{\rm c}^{\prime}(x^{\nu};x - x^{\nu}) \, \geq \, 0,
\epc \forall \, x \, \in \, X.
\end{equation}
In the following, we restrict $\nu$ so that $\zeta_{\rm c}(x^{\nu}) > 0$.  On one hand, we have
\[
\zeta_{\rm c}^{\prime}(x^{\nu};x - x^{\nu}) \, = \, \phi_{\rm c}^{\prime}(x^{\nu};x-x^{\nu}) - ( \, g^{\nu} \, )^T( \, x - x^{\nu} \, ),
\]
where $g^{\nu} \in \partial \varphi_{\rm c}(x^{\nu})$ that depends on the vector $x$.  With $x$ fixed,
the sequence of subgradients $\{ g^{\nu} \}$ has an accumulation
point $g^{\infty}$ which belongs to $\partial \varphi_{\rm c}(x^{\infty})$.  Without loss of generality, we may assume
that $g^{\infty}$ is the limit of the sequence $\{ g^{\nu} \}$.  Consequently, we deduce that
\[
\displaystyle{
\limsup_{\nu \to \infty}
} \, \zeta_{\rm c}^{\prime}(x^{\nu};x - x^{\nu}) \, \leq \, \phi_{\rm c}^{\prime}(x^{\infty};x-x^{\infty}) - ( \, g^{\infty} \, )^T( \, x - x^{\infty} \, ).
\]
On the other hand, by the boundedness of $X$ and the global Lipschitz continuity of $\zeta$, it follows that
$\zeta^{\prime}(x^{\nu};x - x^{\nu})$ is bounded.  Hence diving by $\rho_{\nu}$
in (\ref{eq:stationary nu}), we deduce
\[
\phi_{\rm c}^{\prime}(x^{\infty};x-x^{\infty}) - ( \, g^{\infty} \, )^T( \, x - x^{\infty} \, ) \, \geq \, 0,
\]
where the limit $g^{\infty}$ depends on $x$.  Thus, we have proved that
\[
\phi_{\rm c}^{\prime}(x^{\infty};x-x^{\infty}) \, \geq \, \displaystyle{
\min_{g \in \partial \varphi_{\rm c}(x^{\infty})}
} \, g^T ( \, x - x^{\infty} \, ), \epc \forall \, x \, \in \, X.
\]
Hence, if $\varphi_{\rm c}$ is strictly differentiable at $x^{\infty}$, the above inequality yields the
d-stationarity of the dc constraint function $\zeta_{\rm c}$ on $X$.

\gap

(b) If $\zeta_{\rm c}(x^{\nu}) < 0$, then the stationarity condition of $x^{\nu}$ is
\[
\zeta^{\prime}(x^{\nu};x - x^{\nu}) \, \geq \, 0, \epc \forall \, x \, \in \, X.
\]
Passing to the limit $\nu \to \infty$ for these $\nu$'s easily yields $\zeta^{\prime}(x^{\infty}; x - x^{\infty}) \geq 0$
for all $x \in X$, as desired.

\gap

(c) Suppose $\zeta_{\rm c}(x^{\nu}) = 0$ for all but finitely many $\nu$'s.  The stationarity condition of $x^{\nu}$ is
\begin{equation} \label{eq:stationary nu zero}
\zeta^{\prime}(x^{\nu};x - x^{\nu}) + \rho_{\nu} \, \max( \, 0, \, \zeta_{\rm c}^{\prime}(x^{\nu};x - x^{\nu}) \, ) \, \geq \, 0,
\epc \forall \, x \, \in \, X,
\end{equation}
from which we want to show: $\zeta^{\prime}(x^{\infty};x - x^{\infty}) \geq 0$ for all
$x \in \wh{Y}^j(x^{\infty})$, where
\[
\wh{Y}^j(x^{\infty}) \, \triangleq \, \left\{ \, x \, \in \, X \, \mid \, \phi_{\rm c}(x) \, \leq \,
\psi_{{\rm c},j}(x^{\infty}) + \nabla \psi_{{\rm c},j}(x^{\infty})^T( x - x^{\infty} ) \, \right\},
\]
with $j$ being the single element of $\in {\cal M}_{\rm c}(x^{\infty})$.  Thus ${\cal M}_{\rm c}(x^{\nu}) = \{ j \}$
also, for all $\nu$ sufficiently large.  Hence,
\[
\zeta_{\rm c}^{\prime}(x^{\nu};x - x^{\nu}) \, = \, \phi_{\rm c}^{\prime}(x^{\nu};x - x^{\nu}) -
\nabla \psi_{\rm c}(x^{\nu})^T ( \, x - x^{\nu} \, ).
\]
Consider a vector $x \in \wh{Y}^j(x^{\infty})$
that satisfies the inequality therein strictly.  Since $\psi_{{\rm c},j}(x^{\infty}) = \varphi_{\rm c}(x^{\infty})
= \phi_{\rm c}(x^{\infty})$, we deduce
\[ \begin{array}{lll}
\displaystyle{
\limsup_{\nu \to \infty}
} \, \phi_{\rm c}^{\, \prime}(x^{\nu};x - x^{\nu}) & \leq & \phi_{\rm c}^{\prime}(x^{\infty};x - x^{\infty}) \\ [5pt]
& \leq & \phi_{\rm c}(x) - \phi_{\rm c}(x^{\infty}) \, <  \, \nabla \psi_{{\rm c},j}(x^{\infty})^T( x - x^{\infty} )
\, = \, \displaystyle{
\lim_{\nu \to \infty}
} \, \nabla \psi_{{\rm c},j}(x^{\nu})^T( x - x^{\nu} ).
\end{array}
\]
It follows that for all $\nu$ sufficiently large,
\[
\phi_{\rm c}^{\, \prime}(x^{\nu};x - x^{\nu}) \, < \, \nabla \psi_{{\rm c},j}(x^{\nu})^T( x - x^{\nu} )
\]
Hence
\[
0 \, \leq \, \zeta^{\prime}(x^{\nu};x - x^{\nu}) + \rho_{\nu} \, \max( \, 0, \, \zeta_{\rm c}^{\prime}(x^{\nu};x - x^{\nu}) \, )
\, = \, \zeta^{\prime}(x^{\nu};x - x^{\nu}).
\]
If $x \in \wh{Y}^j(x^{\infty})$ is such that $\phi_{\rm c}(x) =
\psi_{{\rm c},j}(x^{\infty}) + \nabla \psi_{{\rm c},j}(x^{\infty})^T( x - x^{\infty} )$,
then the vector $x^{\tau} \triangleq \bar{x}^j + \tau ( x - \bar{x}^j )$, where
$\bar{x}^j$ is the Slater point under the CQ, i.e., $\bar{x}^j \in X$ satisfies the dc constraint
strictly, remains a Slater point of $\wh{Y}^j$ for all $\tau \in [ 0, 1 )$.  By the above proof, we have
$\zeta^{\prime}(x^{\nu};x^{\tau} - x^{\nu}) \geq 0$.  Letting $\tau \uparrow 1$ completes the proof.  \hfill $\Box$

\gap

{\bf Remarks.} The assumption that ${\cal M}_{\rm c}(x^{\infty})$ is a singleton, or equivalently that
$\varphi_{\rm c}$ is strictly differentiable at $x^{\infty}$, is a pointwise goodness of the
dc function $\zeta_{\rm c}(x) = \phi_{\rm c}(x) - \varphi_{\rm c}(x)$ at $x^{\infty}$.  In spite of the
differentiability of the function $\varphi_{\rm c}$ at $x^{\infty}$, the difference function $\zeta_{\rm c}$ remains not
necessarily so.  This is another instance where the class of good dc functions offers an advantage over the class
of not-good dc functions.

\gap

Another noteworthy remark is that while the assumption $\zeta_{\rm c}(x^{\nu}) = 0$ renders $x^{\nu}$ feasible to the
problem (\ref{eq:dc program dc constraints}), this vector is obtained as a limit point of a presumably infinite
process when the penalized subproblem is solved, and is thus generally not readily available in practical computation.
From this perspective, Proposition~\ref{pr:penalized dc} should be considered
a conceptual result in that it offers insights into the asymptotic property of the penalization approach
for solving a dc constrained dc program without assuming feasibility. How this result can be turned into a constructive
approach for use in practice in solving such a problem requires further investigation.  \hfill $\Box$

\section{Algorithmic Extension: II} \label{sec:extension II}

In this section, we discuss how we can develop a distributed algorithm for solving the following extended dc program:
\begin{equation} \label{eq:dc program sum obj}
\displaystyle{
\operatornamewithlimits{\mbox{minimize}}_{x \in X}
} \ \zeta(x) \, \triangleq \, \phi(x) - \displaystyle{
\sum_{i=1}^I
} \, \varphi_i(x), \epc \mbox{with each} \
\varphi_i(x) \, \triangleq \, \displaystyle{
\operatornamewithlimits{\max}_{1 \leq k \leq \ell_i}
} \, \psi_{i,k}(x)
\end{equation}
where $\phi$ and each $\psi_{i,k}$ are
convex functions defined on $\Omega$ with each $\psi_{i,k}$ being C$^1$ on $\Omega$.  The goal is
to exploit the sum structure in the objective function so that each summand can be treated
separately from the others.  One motivation of this consideration arises from a multi-agent optimization context wherein
each agent $i$ has a private performance function $\varphi_i(x)$ and it is desirable to be able to
implement an algorithm requiring minimal communication among the agents.  Two challenges of this goal
is the dc and nondifferentiable features of the overall objective function and the
coupling of variables in each summand.

\gap

Before presenting the distributed algorithm, we mention that while it is possible to apply the basic
Algorithm~I to the problem (\ref{eq:dc program sum obj}), the sum structure makes a straightforward
application of this centralized algorithm rather laborious when there are many summands.  In this
case, it may be necessary to solve many subproblems at each iteration that are derived by selecting
the functions $\psi_{i,j}$ for all
$j \in \left\{ k \mid \psi_{i,k}(x^{\nu}) \geq \varphi_i(x^{\nu}) - \varepsilon \right\}$
and all $i = 1, \cdots, I$.   To give an example, consider the case $I = n$ and each $\ell_i = 2$.
In this case, the number of subproblems to be solved in each iteration could be exponential in $n$.
Randomization could help in this regard by not exhausting such a selection
per iteration; yet the resulting algorithm remains a centralized scheme that does not take advantage of
the sum structure for possible parallel processing.  To see how the probabilistic
approach can be applied, we define the tuple $t \triangleq \left( k_i \right)_{i=1}^I$ and let
${\cal K} \triangleq \displaystyle{
\prod_{i=1}^I
} \, \{ 1, \cdots, \ell_i \}$.  For each such tuple $t$, let $\psi_t(x) \triangleq \displaystyle{
\sum_{i=1}^I
} \, \psi_{i,k_i}(x)$.  It is easy to see that
\[
\zeta(x) \, = \, \phi(x) - \varphi(x), \epc \mbox{where }
\varphi(x) \, \triangleq \, \displaystyle{
\max_{t \in {\cal K}}
} \, \psi_t(x).
\]
The total number of elements in ${\cal K}$ is $\displaystyle{
\prod_{i=1}^I
} \, \ell_i$, which could be very large.  [For instance, if each $\ell_i = 2$ and $I = n$, then the product
of these $\ell_i$'s is equal to $2^n$, which is exponential in the dimension of the variable $x$.]  In this case,
the randomized version of the algorithm becomes useful.  For $\varepsilon > 0$ and $i = 1, \cdots, I$, let
\[
{\cal M}_i(x) \, \triangleq \, \{ \, k \, \mid \, \psi_{i,k}(x) \, = \, \varphi_i(x) \, \}
\epc \mbox{and} \epc
{\cal M}_{i,\varepsilon}(x) \, \triangleq \, \{ \, k \, \mid \, \psi_{i,k}(x) \, \geq \, \varphi_i(x) - \varepsilon \, \}.
\]
% As in the previous algorithm, we prescribe an index $\iota_i(x) \in {\cal M}_i(x)$ corresponding to a given vector $x$.
At each iteration $\nu$, given an iterate $x^{\nu}$, we
% let $\mathbf{\bar{k}}_{\nu} \triangleq \left( \bar{k}_{\nu,i} \right)_{i=1}^I$,
% where $\bar{k}_{\nu,i} \triangleq \iota_i(x^{\nu})$.
% In addition to this prescribed tuple, we
select a random tuple $t \triangleq \left( k_i \right)_{i=1}^I \in {\cal K}$
such that $k_i \in {\cal M}_{i,\varepsilon}(x^{\nu})$ for every $i = 1, \cdots, I$.  We then solve the following
strongly convex subproblem,
\begin{equation} \label{eq:sum subproblem}
% \begin{array}{ll}
\displaystyle{
\operatornamewithlimits{\mbox{argmin}}_{x \in X}
} \ \left[ \, \phi(x) + \thalf \, \| \, x - x^{\nu} \, \|^2 - \displaystyle{
\sum_{i=1}^I
} \, \left(  \psi_{k_i}(x^{\nu}) + \nabla \psi_{k_i}(x^{\nu})^T( x - x^{\nu} ) \, \right\} \, \right].
\end{equation}
Although the randomized selection of the tuple $t$ avoids the enumeration of a possibly large number of elements
of the set ${\cal K}$ and significantly reduces the number of subproblems to be solved at each iteration, the global
resolution of (\ref{eq:sum subproblem}) remains a centralized task.
While it may be possible to simplify this task under some structural assumptions on the set $X$ and differentiability
properties of the functions $\phi(x)$ and $\varphi(x)$ (see \cite{Alvarado14,ASPang14}), we present below a distributed penalty
approach that is by itself a novel idea for computing d-stationary points of dc programs of the kind (\ref{eq:dc program sum obj})
and requires no such additional structures.  Variations of this approach can be applied to other separable forms of a dc program
(e.g., when $\phi(x)$ is also a sum of agents' functions or a sum of a differentiable and a non-differentiable function).
In what follows, we restrict our discussion to (\ref{eq:dc program sum obj})
where a sum structure is present only in the concave term of the objective.  This distributed approach recognizes
the sum structure $\displaystyle{
\sum_{i=1}^I
} \, \displaystyle{
\operatornamewithlimits{\max}_{1 \leq k \leq \ell_i}
} \, \psi_{i,k}(x)$ and solves subproblems that naturally decomposes according to the latter structure; each decomposed subproblem
can be solved in parallel per individual summand.

\subsection{A penalty approach}

The penalty approach for computing a d-stationary solution to the problem (\ref{eq:dc program sum obj})
consists of two main iterative steps implemented by a sequence of outer iterations each in turn composed of a sequence
of inner iterations.  Each outer iteration is based on the simple observation that
the problem (\ref{eq:dc program sum obj})
is clearly equivalent to the following one where the single variable $x$ is duplicated
$I$ times with the addition of the constraints: $z^i = x$.  This results in a reformulated problem
with $I + 1$ variables:
\begin{equation} \label{eq:dc program sum obj reform}
\begin{array}{ll}
\displaystyle{
\operatornamewithlimits{\mbox{minimize}}_{x, \, z^i \, \in \, X}
} & \phi(x) - \displaystyle{
\sum_{i=1}^I
} \, \varphi_i(z^i) \\ [0.2in]
% \epc \mbox{with each} \ \varphi_i(z^i) \, \triangleq \, \displaystyle{
% \operatornamewithlimits{\max}_{1 \leq k \leq \ell_i}
% } & \psi_{i,k}(x^i) \\ [0.3in]
\mbox{subject to} & z^i \, = \, x, \epc i \, = \, 1, \cdots, I.
\end{array}
\end{equation}
We next penalize the duplication constraints by replacing them with a sum-of-squares term in the objective
using a penalty scalar $\rho > 0$:
\begin{equation} \label{eq:dc program sum obj penalty}
\displaystyle{
\operatornamewithlimits{\mbox{minimize}}_{x, \, z^i \, \in \, X}
} \ \theta_{\rho}(x,z) \, \triangleq \, \phi(x) - \displaystyle{
\sum_{i=1}^I
} \, \varphi_i(z^i) + \displaystyle{
\frac{\rho}{2}
} \, \displaystyle{
\sum_{i=1}^I
} \, \left\| \, z^i - x \, \right\|^2; \mbox{ where }
z \, \triangleq \, \left( \, z^i \, \right)_{i=1}^I.
% \epc \mbox{with each} \
% \varphi_i(x) \, \triangleq \, \displaystyle{
% \operatornamewithlimits{\max}_{1 \leq k \leq \ell_i}
% } \ \psi_{i,k}(z^i)
\end{equation}
The outer iterations consist of solving the problem (\ref{eq:dc program sum obj penalty}) for an increasing
sequence of positive scalars $\{ \rho_{\nu} \}$ tending to $\infty$.
This is accomplished by applying the basic Algorithm I or its randomized version to the following problem:
\begin{equation} \label{eq:dc program sum obj penaltyII}
\displaystyle{
\operatornamewithlimits{\mbox{minimize}}_{x, \, z^i \, \in \, X}
} \ \theta_{\rho_{\nu}}(x,z) \, \triangleq \, \phi(x) - \displaystyle{
\sum_{i=1}^I
} \, \varphi_i(z^i) + \displaystyle{
\frac{\rho_{\nu}}{2}
} \, \displaystyle{
\sum_{i=1}^I
} \, \left\| \, z^i - x \, \right\|^2
\end{equation}
for each $\nu$, yielding a sequence of d-stationary points $\left\{ x^{\nu}, ( z^{\nu,i} )_{i=1}^I \right\}_{\nu=1}^{\infty}$.
Thus, for all $x$ and $z^i$ in $X$,
\begin{equation} \label{eq:penal stationary}
\begin{array}{l}
\phi^{\, \prime}(x^{\nu};x - x^{\nu}) - \displaystyle{
\sum_{i=1}^I
} \, \displaystyle{
\max_{k \in {\cal M}_i(z^{\nu,i})}
} \, \nabla \psi_{i,k}(z^{\nu,i})^T ( \, z^i - z^{\nu,i} \, ) + \\ [0.1in]
\epc \rho_{\nu} \, \displaystyle{
\sum_{i=1}^I
} \, \left[ \, ( \, z^{\nu,i} - x^{\nu} \, )^T( \, z^i - z^{\nu,i} \, ) + ( \, x^{\nu} - z^{\nu,i} \, )^T( \, x - x^{\nu} \, ) \, \right]
\, \geq \, 0.
\end{array}
\end{equation}
Before describing the inner iterations to generate such stationary solutions of (\ref{eq:dc program sum obj penaltyII}),
we first establish the desired limiting property of such solutions; namely, every accumulation point
$\left( x^{\infty}, (z^{\infty,i})_{i=1}^I \right)$ of
$\left\{ x^{\nu}, ( z^{\nu,i} )_{i=1}^I \right\}_{\nu=1}^{\infty}$ must satisfy $z^{\infty,i} = x^{\infty}$ for all $i = 1, \cdots,I$;
thus we recover the feasibility condition of (\ref{eq:dc program sum obj reform}).
% More importantly, we show that $x^{\infty}$ is a d-stationary solution of (\ref{eq:dc program sum obj}).

\gap

{\bf Convergence of penalization}.
Throughout the following analysis, we assume that % $X$ is bounded,
each $\| \nabla \psi_{i,k} \|$ is bounded on $X$.
% and $\phi$ is Lipschitz continuous on $X$ with a Lipschitz constant $L > 0$ so that
% $| \phi^{\, \prime}(x;d) | \leq L \| d \|$ for all $x \in X$ and all $d \in \mathbb{R}^n$.
% The stationarity condition (\ref{eq:penal strong stationary}) implies
% \[ % begin{equation} \label{eq:variational nu}
% \phi^{\, \prime}(x^{\nu};x - x^{\nu}) - \displaystyle{
% \sum_{i=1}^I
% } \, \displaystyle{
% \max_{k \in {\cal M}_i(z^{\nu,i})}
% } \, \nabla \psi_{i,k}(z^{\nu,i})^T ( \, z^i - z^{\nu,i} \, ) + \displaystyle{
% \frac{\rho_{\nu}}{2}
% } \, \displaystyle{
% \sum_{i=1}^I
% } \, \left[ \, \left\| \, z^i - x \, \right\|^2 - \| \, z^{\nu,i} - x^{\nu} \, \|^2 \, \right] \, \geq \, 0.
% \]
By letting $\eta \triangleq \displaystyle{
\max_{1 \leq i \leq I}
} \, \displaystyle{
\max_{1 \leq k \leq \ell_i}
} \, \displaystyle{
\max_{x \in X}
} \, \| \nabla \psi_{i,k}(x) \|$ and $z^i = x = x^{\nu}$ for all $i$, we deduce from (\ref{eq:penal stationary}),
\[
0 \, \leq \, \eta \, \displaystyle{
\sum_{i=1}^I
} \, \| \, x^{\nu} - z^{\nu,i} \, \| - \displaystyle{
\frac{\rho_{\nu}}{2}
} \, \displaystyle{
\sum_{i=1}^I
} \, \| \, z^{\nu,i} - x^{\nu} \, \|^2,
\]
% With the choice of $z^i = x$ for all $i$, it follows that
% \[
% \displaystyle{
% \frac{\rho_{\nu}}{2}
% } \, \displaystyle{
% \sum_{i=1}^I
% } \, \| z^{\nu,i} - x^{\nu} \|^2 \, \leq \, L \, \| \, x - x^{\nu} \, \| + \eta \, \displaystyle{
% \sum_{i=1}^I
% } \, \| \, x - z^{\nu,i} \, \| .
% \]
% Since the right-hand side is bounded, so is the left-hand side.
% Consequently,
which easily implies that $\displaystyle{
\lim_{\nu \to \infty}
} \, \| z^{\nu,i} - x^{\nu} \| = 0$ for all $i$.  Hence, if $x^{\infty}$ is the limit of
a convergent subsequence $\{ x^{\nu} \}_{\nu \in {\cal N}}$,
which must exist by the boundedness of $X$, then $\displaystyle{
\lim_{\nu (\in {\cal N}) \to \infty}
} \, z^{\nu,i} = x^{\infty}$ for all $i$.  With $z^i = x$ for every $i$, (\ref{eq:penal stationary})
also implies
\[ % begin{equation} \label{eq:variational nu}
\phi^{\, \prime}(x^{\nu};x - x^{\nu}) \, \geq \, \displaystyle{
\sum_{i=1}^I
} \, \displaystyle{
\max_{k \in {\cal M}_i(z^{\nu,i})}
} \, \nabla \psi_{i,k}(z^{\nu,i})^T ( \, x - z^{\nu,i} \, ).
\]
Since ${\cal M}_i(z^{\nu,i}) \subseteq {\cal M}_i(x^{\infty})$ for all $\nu \in {\cal N}$
sufficiently large, we deduce that for some nonnegative scalars $\left\{ \lambda_{i,k}^{\nu} \right\}_{k \in {\cal M}_i(x^{\infty})}$,
satisfying $\displaystyle{
\sum_{k \in {\cal M}_i(x^{\infty})}
} \, \lambda_{i,k}^{\nu} = 1$ and possibly dependent on $x$,
\[
\phi^{\, \prime}(x^{\nu};x - x^{\nu}) \, \geq \, \displaystyle{
\sum_{i=1}^I
} \, \displaystyle{
\sum_{k \in {\cal M}(x^{\infty})}
} \, \lambda_{i,k}^{\nu} \, \nabla \psi_{i,k}(z^{\nu,i})^T ( \, x - z^{\nu,i} \, ).
\]
For $x$ fixed, we may assume, without loss of generality, that for each pair $(i,k)$, the sequence of scalars
$\left\{ \lambda_{i,k}^{\nu} \right\}_{\nu \in {\cal N}}$ converges to $\lambda_{i,k}^{\infty}$,
which must be nonnegative and satisfies: $\displaystyle{
\sum_{k \in {\cal M}_i(x^{\infty})}
} \, \lambda_{i,k}^{\infty} = 1$.
By a known limiting property of the directional derivatives of convex functions \cite[Theorem~24.5]{Rockafellar70}, we have
\[
\phi^{\, \prime}(x^{\infty};x - x^{\infty}) \, \geq \,
\displaystyle{
\limsup_{\nu (\in {\cal N}) \to \infty}
} \, \phi^{\, \prime}(x^{\nu};x - x^{\nu}).
\]
Hence, % for all $x \in X$
\[
\phi^{\, \prime}(x^{\infty};x - x^{\infty}) \, \geq \,  \displaystyle{
\sum_{i=1}^I
} \, \displaystyle{
\sum_{k \in {\cal M}_i(x^{\infty})}
} \, \lambda_{i,k}^{\infty} \, \nabla \psi_{i,k}(x^{\infty})^T ( \, x - x^{\infty} \, ).
\]
Since $\displaystyle{
\sum_{k \in {\cal M}_i(x^{\infty})}
} \, \lambda_{i,k}^{\infty} \nabla \psi_{i,k}(x^{\infty}) \in \partial \varphi_i(x^{\infty})$, we deduce that
\[
\phi^{\, \prime}(x^{\infty};x - x^{\infty}) \, \geq \,  \displaystyle{
\sum_{i=1}^I
} \, \displaystyle{
\operatornamewithlimits{
\min_{g^i \in \partial \varphi_i(x^{\infty})}}
} \, ( \, g^i \, )^T ( \, x - x^{\infty} \, ), \epc \forall \, x \, \in \, X.
\]
Hence, if each $\partial \varphi_i(x^{\infty})$ is a singleton, it follows that
$x^{\infty}$ is a d-stationarity solution of (\ref{eq:dc program sum obj}).  We have therefore proved
the next result.

\begin{proposition} \label{pr:convergence distributed} \rm
Suppose that each $\| \nabla \psi_{i,k} \|$ is bounded on $X$.
\begin{description}
\item[\rm (a)] (Recovering feasibility) Every accumulation point
$\left( x^{\infty}, (z^{\infty,i})_{i=1}^I \right)$
of the sequence
$\left\{ x^{\nu}, ( z^{\nu,i} )_{i=1}^I \right\}_{\nu=1}^{\infty}$
of penalized d-stationary points corresponding to a sequence of penalty
parameters $\{ \rho_{\nu} \} \uparrow \infty$ must satisfy $z^{\infty,i} = x^{\infty}$ for all $i = 1, \cdots I$.
\item[\rm (b)] (Achieving stationarity)
Moreover, if each $\varphi_i$ is strictly differentiable at $x^{\infty}$, then $x^{\infty}$
is a d-stationary solution of (\ref{eq:dc program sum obj}).  \hfill  $\Box$
\end{description}
\end{proposition}

{\bf Remark.}  Once again, the goodness of the objective function of (\ref{eq:dc program sum obj}) is
needed to complete the last step of the proof of the above proposition.  \hfill $\Box$

\gap

{\bf A distributed algorithm for (\ref{eq:dc program sum obj penaltyII})}.
Based on Algorithm I, we present in this section a distributed algorithm for computing
a d-stationary solution of each penalized problem (\ref{eq:dc program sum obj penaltyII}).
To do this, we need to take care of one detail of this problem having to do with the non-separability of
the term $\| z^i - x \|^2$ in the objective function.  Namely, we linearize this term at
a base tuple $\left( x^{\nu}, ( z^{\nu,i} )_{i=1}^I \right)$ as follows:
\begin{equation} \label{eq:nonseparable}
% \begin{array}{lll}
\| \, z^i - x \|^2 \, \approx \, \| \, z^{\nu,i} - x^{\nu} \, \|^2 + 2 \, \left[ \, ( \, x - x^{\nu} \, )^T( \, x^{\nu} - z^{\nu,i} \, )
+ ( \, z^i - z^{\nu,i} \, )^T( \, z^{\nu,i} - x^{\nu} \, ) \, \right]
% & = & \| \, z^i \, \|^2 - 2 x^Tz^i + \| \, x \, \|^2 \\ [5pt]
% & \approx & \| \, z^i \, \|^2 - 2 \,
% \left[ \, ( \, x^{\nu} \, )^Tz^{\nu,i} + ( \, x - x^{\nu} \, )^Tz^{\nu,i} + ( \, x^{\nu} \, )^T ( \, z^i - z^{\nu,i} \, ) \, \right]
% + \| \, x \, \|^2 \\ [5pt]
% & = & \underbrace{( \, z^i - x^{\nu} \, )^T( \, z^i - x^{\nu} \, ) + ( \, x - z^{\nu,i} \, )^T( \, x - z^{\nu,i} \, ) -
% ( \, x^{\nu} - z^{\nu,i} \, )^T( \, x^{\nu} - z^{\nu,i} \, )}_{\mbox{separable and strongly convex
% in $z^i$ and $x$ with $( x^{\nu},z^{\nu,i} )$ given}}.
% \end{array}
\end{equation}
and use this linearization in each step of the algorithm.

\gap

At the beginning of an outer iteration $\nu$ (thus $\rho_{\nu}$ is fixed), starting at a tuple
$\left( x^{\nu,0}, \left( z^{\nu,i,0} \right)_{i=1}^I \right)$ of vectors in $X$, the algorithm generates a sequence
of inner iterates
\begin{equation} \label{eq:sum sequence}
\left\{ \, x^{\nu,\mu},( z^{\nu,i,\mu} )_{i=1}^I \, \right\}_{\mu=0}^{\infty} .
\end{equation}
At each inner iteration $\mu = 1, 2, \cdots$, for every tuple
$t_{\nu,\mu} \triangleq \left( k_{\nu,i,\mu} \right)_{i=1}^I$
consisting of indices $k_{\nu,i,\mu} \in {\cal M}_{i,\varepsilon}(z^{\nu,i,\mu})$ for $i = 1, \cdots, I$,
% where $\varepsilon > \varepsilon^{\, \prime}$,
we solve the strongly convex subprogram:
\begin{equation} \label{eq:subproblem nu}
\begin{array}{l}
\displaystyle{
\operatornamewithlimits{\mbox{minimize}}_{x, \, z^i \, \in \, X}
} \ \left\{ \, \phi(x) + \thalf \, \left( \, \underbrace{\| \, x - x^{\nu,\mu} \, \|^2 + \displaystyle{
\sum_{i=1}^I
} \, \| \, z^i - z^{\nu,i,\mu} \, \|^2}_{\mbox{regularization}} \right) + \right. \\ [0.5in]
\hspace{0.5in} \rho_{\nu} \, \displaystyle{
\sum_{i=1}^I
} \, \left[ \, \underbrace{( \, x - x^{\nu,\mu} \, )^T( \, x^{\nu,\mu} - z^{\nu,i,\mu} \, ) +
( \, z^i - z^{\nu,i\,\mu} \, )^T( \, z^{\nu,i,\mu} - x^{\nu,\mu} \, )}_{\mbox{linear approximation of $\thalf \| z^i - x \|^2$}} \, \right] -
\\ [0.35in]
\hspace{1.5in}  \left. \displaystyle{
\sum_{i=1}^I
} \,
\left( \, \psi_{i,k_{\nu,i,\mu}}(z^{\nu,i,\mu}) + \nabla \psi_{i,k_{\nu,i,\mu}}(z^{\nu,i,\mu})^T( \, z^i - z^{\nu,i,\mu} \, ) \, \right)
\, \right\},
\end{array} \end{equation}
which naturally decomposes into $I + 1$ subproblems:

\gap

$\bullet $ a strongly convex subproblem in the $x$-variable,
\[\displaystyle{
\operatornamewithlimits{\mbox{minimize}}_{x \in X}
} \ \left[ \, \phi(x) + \thalf \, \| \, x - x^{\nu,\mu} \, \|^2 +
\rho_{\nu} \, ( \, x - x^{\nu,\mu} \, )^T \, \displaystyle{
\sum_{i=1}^I
} \, ( \, x^{\nu,\mu} - z^{\nu,i,\mu} \, ) \, \right];
\]
$\bullet $ a problem of the same kind in the $z^i$-variable, for $i = 1, \cdots, I$,
\begin{equation} \label{eq:subproblem zi}
\begin{array}{l}
\displaystyle{
\operatornamewithlimits{\mbox{minimize}}_{z^i \in X}
} \ \left\{ \, \thalf \| \, z^i - z^{\nu,i,\mu} \, \|^2 + \rho_{\nu} \, ( \, z^i - z^{\nu,i\,\mu} \, )^T( \, z^{\nu,i,\mu} - x^{\nu,\mu} \, ) 
- \right. \\ [0.2in]
\hspace{1in} \left. \left[ \,
\psi_{i,k_{\nu,i,\mu}}(z^{\nu,i,\mu}) + \nabla \psi_{i,k_{\nu,i,\mu}}(z^{\nu,i,\mu})^T( \, z^i - z^{\nu,i,\mu} \, )
\, \right] \, \right\}.
\end{array}
\end{equation}
% Writing $z^{\nu,\mu} \triangleq \left( z^{\nu,i,\mu} \right)_{i=1}^I$, we let
% \begin{equation} \label{eq:numu}
% \left\{ \, \wh{x}^{\, \mathbf{k}_{\nu,\mu}}(x^{\nu,\mu},z^{\nu,\mu}),
% \left( \, \wh{z}^{\, i,\mathbf{k}_{\nu,\mu}}(x^{\nu,\mu},z^{\nu,\mu}) \, \right)_{i=1}^I \, \right\}
% \end{equation}
% denote an optimal solution to (\ref{eq:subproblem nu}).
Among the optimal solutions to (\ref{eq:subproblem nu}) for various choices of the index tuples $t_{\nu,\mu}$, one of them leads to
the new iterates $\left( x^{\nu,\mu+1}, ( z^{\nu,i,\mu+1} )_{i=1}^I \right)$.
[This is the deterministic version of the algorithm; we leave the probabilistic version for the reader to complete
and remind the reader that the latter could be more efficient than the former in practice especially when synchronization 
among agents is costly.]  In total, for each tuple
$t$, a total of $I + 1$ strongly convex subprograms are solved at each inner iteration, each of them
can be solved separately
from the others.  The choice of the individual indices $k_{\nu,i,\mu}$ can be carried out in parallel per agent $i$.
Thus, the overall implementation of the algorithm
% for computing a d-stationary solution of the strong kind of each penalized
% program (\ref{eq:dc program sum obj penaltyII})
is totally distributed.  At the completion of the inner iterations according to
a prescribed termination criterion, the penalty parameter $\rho_{\nu}$ is updated and a new sequence of inner iterations is entered.
The convergence of the inner iterations can be proved in a similar way as the basic Algorithm~I and is not repeated.

\gap

The outer-inner scheme for solving (\ref{eq:dc program sum obj}) has its advantage of being implementable distributedly according to
the individual summands.  Ideally, it would be desirable to have a single-loop algorithm wherein the update of the penalty
parameter $\rho$ can be incorporated into the inner iterations.  At this time, we are not able to develop a provably convergent
single-loop algorithm that can be implemented distributedly.

\gap

{\bf Acknowledgement.}  The authors are grateful to two referees for offering many constructive comments that have improved the
presentation of the paper.  Moreover, the first author has benefitted from a fruitful visit to Lorraine University where he had
a very productive discussion with Professors Pham Dinh Tao and Le Thi Hoa An on dc programming in general and this paper in particular.

{\bf Appendix: Proof of the dc-property of $\theta$ in (\ref{eq:general_problem})}.
We first introduce some notational definitions.  For every $i=1,\ldots,I$ and $j=1,\ldots,J$ we let
\begin{equation}\label{eq:rho_def}
\rho_{i,j}^{\max} \, \triangleq \, \underset{\lambda^i \in \Lambda^i}{\text{maximum}} \,\, h_{i,j} (\lambda^i)
\epc \mbox{and} \epc
\rho_{i,j}^{\min} \, \triangleq \, \underset{\lambda^i \in \Lambda^i}{\text{minimum}} \,\, h_{i,j} (\lambda^i),
\end{equation}
which are finite scalars by the compactness of $\Lambda^i$ and the continuity of $h_{i,j}$.
Based on these two extremum values, we rewrite each product $h_{i,j} (\lambda^i) f_{i,j}(x)$
as follows:

\gap

{\bf (i)} if $f_{i,j}(x)$ is convex (cvx):
\[
h_{i,j}(\lambda^i) \, f_{i,j}(x) \, = \,
\rho_{i,j}^{\min} \, f_{i,j}(x) + \left( \, h_{i,j}(\lambda^i) - \rho_{i,j}^{\min} \, \right) \, f_{i,j}(x),
\]
{\bf (ii)} if $f_{i,j}(x)$ is concave (cve):
\[ 		
h_{i,j}(\lambda^i) \, f_{i,j}(x) \, = \,
\rho_{i,j}^{\max} \, f_{i,j}(x) + \left( \, \rho_{i,j}^{\max} - h_{i,j}(\lambda^i) \, \right) \, (-f_{i,j}(x)).
\]
Thus, we can rewrite each $\theta_i(x)$ for $i=1,\ldots,I$ as
\begin{equation}\label{eq:theta_i_rewrite}
\begin{aligned}
	\theta_i(x) &=
	\underset{\lambda^i \in \Lambda^i}{\text{maximum}} \,\sum_{j=1}^{J} h_{i,j} (\lambda^i) f_{i,j}(x)\\ %%%
	& = \sum_{j: f_{i,j} \text{ cvx}} \rho_{i,j}^{\min} \, f_{i,j}(x) +
	      \sum_{j: f_{i,j} \text{ cve}} \rho_{i,j}^{\max} \, f_{i,j}(x)  \\
	 & \quad + \underset{\lambda^i \in \Lambda^i}{\text{maximum}}\left(
	\sum_{j: f_{i,j} \text{ cvx}} \left(
			h_{i,j}(\lambda^i) - \rho_{i,j}^{\min}
		\right)f_{i,j}(x) +
	\sum_{j: f_{i,j} \text{ cve}} \left(
			\rho_{i,j}^{\max} - h_{i,j}(\lambda^i)
		\right)(-f_{i,j}(x))
	\right)             %%			
\end{aligned}		
\end{equation}
It is important to highlight the following facts with regard to the above representation:
% that are key points in deriving an iterative algorithm for solving (\ref{eq:general_problem}):

\gap

$\bullet $ If $\rho_{i,j}^{\min} \leq 0$ for some $j$ such that $f_{i,j}$ is convex then
$\rho_{i,j}^{\min} \, f_{i,j}(x)$ is a {\sl concave} function on $X$.

\gap

$\bullet $ If $\rho_{i,j}^{\max} \geq 0$ for some $j$ such that $f_{i,j}$ is concave then
$\rho_{i,j}^{\max} \, f_{i,j}(x)$ is a {\sl concave} function on $X$.

\gap

$\bullet $ Let
\begin{equation}\label{eq:max_function}
g_i(x) \, \triangleq \, \underset{\lambda^i \in \Lambda^i}{\text{maximum}}
\left[ \, \underbrace{\displaystyle{
\sum_{j: f_{i,j} \text{ cvx}}
} \, \left( \, h_{i,j}(\lambda^i) - \rho_{i,j}^{\min} \right) \, f_{i,j}(x)
+ \displaystyle{
\sum_{j: f_{i,j} \text{ cve}}
} \, \left( \, \rho_{i,j}^{\max} - h_{i,j}(\lambda^i) \, \right) \, (-f_{i,j}(x))}_{\mbox{$\triangleq \,
\varphi_i(x,\lambda^i)$}}
\, \right].
\end{equation}
Since the maximand $\varphi_i(x,\lambda^i)$ is a convex function in $x$ for each $\lambda^i$,
it follows readily that $g_i(x)$ is a convex function on $X$.

\gap

$\bullet $ By combining the above three observations, it is clear that
each $\theta_i(x)$ is a dc-function.  To be explicit, we define some index sets:
\[ \begin{array}{lll}
\mathcal{J}_i^{-,{\rm cvx}} & \triangleq & \left\{ \, j \, \mid \, \rho_{i,j}^{\min} \, \leq \, 0 \text{ and }
f_{i,j} \text{ is cvx } \right\} \\ [5pt]
\mathcal{J}_i^{+,{\rm cvx}} & \triangleq & \left\{ \, j \, \mid \, \rho_{i,j}^{\min} \, \geq \, 0 \text{ and }
f_{i,j} \text{ is cvx } \right\} \\ [5pt]
\mathcal{J}_i^{+,{\rm cve}} & \triangleq & \left\{ j \, \mid \, \rho_{i,j}^{\max} \, \geq \, 0 \text{ and }
f_{i,j} \text{ is cve } \right\} \\ [5pt]
\mathcal{J}_i^{-,{\rm cve}} & \triangleq & \left\{ j \, \mid \, \rho_{i,j}^{\max} \, \leq \, 0 \text{ and }
f_{i,j} \text{ is cve } \right\}.
\end{array}
\]
We can then write
\[ \begin{array}{l}
\theta_i(x) \, = \, \displaystyle{
\sum_{j: f_{i,j} \text{ cvx}}
} \, \rho_{i,j}^{\min} \, f_{i,j}(x) + \displaystyle{
\sum_{j: f_{i,j} \text{ cve}}
} \, \rho_{i,j}^{\max} \, f_{i,j}(x) + g_i(x) \\ [0.3in]
= \, \left[ \, \underbrace{\displaystyle{
\sum_{j \in \mathcal{J}_i^{-,{\rm cvx}}}
} \, \rho_{i,j}^{\min} \, f_{i,j}(x) + \displaystyle{
\sum_{j \in \mathcal{J}_i^{+,{\rm cve}}}
} \, \rho_{i,j}^{\max} \, f_{i,j}(x)}_{\mbox{$\triangleq \, u_i(x)$}} \, \right] +
\left[ \, \underbrace{\displaystyle{
\sum_{j \in \mathcal{J}_i^{+,{\rm cvx}}}
} \, \rho_{i,j}^{\min} \, f_{i,j}(x) + \displaystyle{
\sum_{j \in \mathcal{J}_i^{-,{\rm cve}}}
} \, \rho_{i,j}^{\max} \, f_{i,j}(x) + g_i(x)}_{\mbox{$\triangleq \, v_i(x)$}} \, \right]
\end{array} \]
where the function $u_i(x)$ is concave and differentiable while $v_i(x)$ is convex and non-differentiable.
Hence $\theta_i(x)$ is a dc function; thus so is $\theta(x) = \displaystyle{
\sum_{i=1}^I
} \, \theta_i(x) = u(x) + v(x)$, where
$u(x) \triangleq \displaystyle{
\sum_{i=1}^I
} \, u_i(x)$ is concave and differentiable and $v(x) \triangleq \displaystyle{
\sum_{i=1}^I
} \, v_i(x)$ is convex and nondifferentiable. Consequently, (\ref{eq:general_problem}) is a nondifferentiable dc
program.  The noteworthy point is that in the representation $\theta(x) = u(x) + v(x)$, the concave summand
$u(x)$ is differentiable whereas the convex summand $v(x)$ is not; thus $\theta$ is not a good dc function.
% it is important to recall that
% the problem (\ref{eq:general_problem}) is the maximization of $\theta(x)$.  It turns out that this feature of $\theta$
% renders the maximization problem (\ref{eq:general_problem}) not a ``good dc program'' in a sense to be made precise later.

\end{document}